\documentclass[reqno]{amsart}
\usepackage{amsfonts}
\usepackage{amssymb}
\usepackage{amsthm}
\usepackage{upref}
\usepackage{enumerate}

\makeatletter
\def\LaTeX{\leavevmode L\raise.42ex
    \hbox{\kern-.3em\size{\sf@size}{0pt}\selectfont A}\kern-.15em\TeX}
 \makeatother
 
\DeclareMathOperator{\clos}{clos}
\DeclareMathOperator{\Det}{Det}

\numberwithin{equation}{section}

\newtheorem{lemma}{Lemma}[section]
\newtheorem{theorem}[lemma]{Theorem} 
\newtheorem{corollary}[lemma]{Corollary}

\theoremstyle{definition}

\newtheorem{example}[lemma]{Example}

\newtheorem{remark}[lemma]{Remark}

 \newcommand{\supp}{\operatorname{supp}}
  \newcommand{\e}{\eqref}

\newcommand{\q}{\quad}

\newcommand{\ov}{\overline}

\newcommand{\z}{\zeta}

\renewcommand{\d}{\delta}

\renewcommand\Im{\operatorname{Im}}
\renewcommand\Re{\operatorname{Re}}

\newenvironment{pf}{\begin{proof}}{\end{proof}}

\def\qqq{\mathrel{\subset\mkern-15mu\lower.38ex\hbox{${\scriptscriptstyle\rightarrow}$}}}

\let\cal\mathcal

\let\Bbb\mathbb
 
\begin{document}
\title 
{Semiclassical asymptotic  behavior of orthogonal polynomials}
\author{ D. R. Yafaev  }
\address{   Univ  Rennes, CNRS, IRMAR-UMR 6625, F-35000
    Rennes, France and SPGU, Univ. Nab. 7/9, Saint Petersburg, 199034 Russia}
\email{yafaev@univ-rennes1.fr}
\subjclass[2000]{33C45, 39A70,  47A40, 47B39}
 
 \keywords {Jacobi matrices,    long-range   perturbations, difference equations, 
 orthogonal polynomials, asymptotics for large numbers.  }

\thanks {Supported by  project   Russian Science Foundation   17-11-01126}

\begin{abstract}
Our goal  is to find  asymptotic formulas for orthonormal polynomials $P_{n}(z)$  with the recurrence coefficients slowly stabilizing as $n\to\infty$. To that end, we develop scattering theory  of    Jacobi operators with long-range coefficients and study the corresponding second order difference equation.  We introduce the Jost solutions $f_{n}(z)$ of this equation by a condition for $n\to\infty$ and suggest an Ansatz for  them  playing the role of the semiclassical Liouville-Green Ansatz for the corresponding solutions of the   Schr\"odinger equation.  This allows us to study Jacobi operators and their eigenfunctions $P_{n}(z)$ by traditional methods of spectral theory developed for differential equations. In particular, we express all coefficients in asymptotic formulas for $P_{n}(z)$ as $\to\infty$ in terms of the Wronskian of the solutions $\{P_{n} (z)\}$ and $\{f_{n} (z)\}$.
   \end{abstract}
   
  \thispagestyle{empty}

\maketitle

\section{Introduction}

\subsection{Jacobi and orthogonal polynomials}

  As is well known, the theories of Jacobi operators given by three-diagonal matrices
  \begin{equation}
J= 
\begin{pmatrix}
 b_{0}&a_{0}& 0&0&0&\cdots \\
 a_{0}&b_{1}&a_{1}&0&0&\cdots \\
  0&a_{1}&b_{2}&a_{2}&0&\cdots \\
  0&0&a_{2}&b_{3}&a_{3}&\cdots \\
  \vdots&\vdots&\vdots&\ddots&\ddots&\ddots
\end{pmatrix} 
\label{eq:ZP+}\end{equation}
 in the space  $\ell^2 ({\Bbb Z}_{+})$  and of differential Schr\"odinger operators $H=D p (x)D +q(x)$ (with, for example, the boundary condition $u(0)=0$)  in the space  $L^2 ({\Bbb R}_{+})$ are  to a large extent similar.   For Jacobi operators, $n\in {\Bbb Z}_{+}$ plays the role of $x\in {\Bbb R}_{+}$ and the coefficients $a_{n}$, $b_{n}$, play the roles of the functions $p(x)$, $q(x)$, respectively.
  In our opinion,  a consistent  analogy between Jacobi and Schr\"odinger operators  sheds a new light on  some aspects of the orthogonal polynomials theory. 
Of course this point of view is not new; for example, it was   advocated long ago by   
    K.~M.~Case in   \cite{Case}.

    In this paper, the sequences $a_{n}>0$ and $b_{n}=\bar{b}_{n}$ in \e{eq:ZP+}  are    assumed to be   bounded, so that $J$ is a bounded self-adjoint operator in   the space $\ell^2 ({\Bbb Z}_{+})$.
    Its spectral family will be denoted    $E(\lambda)$. 
 The spectrum of $J$ is simple with $e_{0} = (1,0,0,\ldots)^\top$ being a generating vector. It is natural to define   the   spectral measure of $J$ by the relation $d\rho(\lambda)=d(E(\lambda)e_{0}, e_{0})$.

Orthogonal polynomials $P_{n } (z)$  associated with the
Jacobi matrix \e{eq:ZP+} are defined  by the recurrence relation
\begin{equation}
 a_{n-1} P_{n-1} (z) +b_{n} P_{n } (z) + a_{n} P_{n+1} (z)= z P_n (z), \q n\in{\Bbb Z}_{+}, 
\label{eq:pol}\end{equation} 
and the boundary conditions $P_{-1 } (z) =0$, $P_0 (z) =1$.    
 Obviously, $P_{n } (z)$ is a polynomial of degree $n$  and the vector $P(z)=\{ P_{n} (z)\}_{n=-1}^\infty $ formally satisfies the equation $JP( z )= z P(z)$, that is, it is an ``eigenvector"  of the operator $H$.
 The  polynomials  $P_{n}(\lambda)$ are orthogonal and normalized  in the space $L^2 ({\Bbb R};d\rho)$: 
      \[
\int_{-\infty}^\infty P_{n}(\lambda) P_{m}(\lambda) d\rho(\lambda) =\d_{n,m};
\]
as usual, $\d_{n,n}=1$ and $\d_{n,m}=0$ for $n\neq m$. 
 Alternatively,   given the probability measure $d\rho(\lambda)$, the polynomials $P_0 (\lambda),P_1 (\lambda),\ldots,P_{n} (\lambda),\ldots$ can be obtained by the Gram-Schmidt orthonormalization of the monomials $1,\lambda,\ldots,\lambda^n,\ldots$ in the space $L^2({\Bbb R}_{+}; d\rho)$; one also has to additionally require that $P_{n} (\lambda)=k_{n}(\lambda^n+ r_{n}\lambda^{n-1} +\cdots)$ with $k_{n} >0$.   The coefficients $a_{n}, b_{n}$ can be recovered by the formulas  $a_{n}= k_{n}  / k_{n+1}$, $b_{n}= r_{n}
 - r_{n+1}$.

     The operator \e{eq:ZP+} with the coefficients $a_{n}=1/2$, $b_{n}=0$ is known as the  ``free"  discrete Schr\"odinger operator. This operator, denoted $J_{0}$,   plays the  role of the  differential operator $D^2$ in the space $L^2 ({\Bbb R}_{+})$.           The   operator $J_{0}$  can   be diagonalized explicitly.
Its spectrum   is  absolutely continuous and coincides with the interval $[-1,1]$. 
The eigenfunctions of $J_{0}$  are   normalized  Chebyshev polynomials   $P_{n}  (\lambda)$  of the second kind, and the corresponding spectral measure $d\rho_{0} (\lambda)= d(E_{0} (\lambda) e_{0}, e_{0})$ is given by the formula
 \begin{equation}
d\rho_{0} (\lambda)= 2 \pi^{-1}  \sqrt{1-\lambda^2 }d\lambda,\q \lambda\in (-1,1). 
\label{eq:fr}\end{equation} 
Obviously, $ J= J_{0}+V $
where the ``perturbation" $V  $ is given by  equality \e{eq:ZP+} with $a_{n}$ replaced by $\alpha_{n}=a_{n}-1/2$.

We refer to the books \cite{AKH, Ism,Sz} and surveys \cite{Lub, Tot} for   general information about orthogonal polynomials.

\subsection{Statement of the problem}

 Our goal is to study an asymptotic behavior of  the polynomials $P_n (z)$ as $n\to\infty$. Of course, one has to distinguish the cases of $z$ in the spectrum $\sigma (J)$ of the Jacobi operator $J$ and of $z \not \in\sigma (J)$. Asymptotic properties of $P_n (z)$ can be deduced either from the coefficients $a_{n}$, $b_{n}$ of the  operator $J$ or from its spectral measure $d\rho (\lambda)$:
 
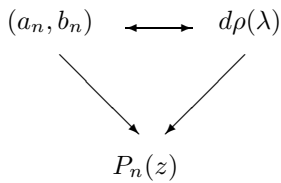
 \begin{picture}(150,80)
 \put(100,55){$(a_{n} , b_{n})$}
\put(180,55){$d\rho(\lambda)$}
\put(170,55){\vector(-1,0){25}}
\put(145,55){\vector(1,0){25}}
\put(120,45){\vector(1,-1){30}}
\put(190,45){\vector(-1,-1){30}}
\put(140, 0){$P_{n} (z)$}
\end{picture}

 \bigskip  \bigskip

Asymptotic formulas were very well known  (see, e.g., the book \cite{BE}) for the classical, that is, Jacobi, Laguerre and Hermite,  polynomials, but the first general result is probably due to S.~Bernstein 
(see his pioneering paper \cite{Bern} or  Theorem~12.1.4 in the G.~Szeg\H{o} book~\cite{Sz}). These results were stated in terms of the measure $d\rho (\lambda)$.
It was required      that $\supp\rho\subset [-1,1]$, the measure is absolutely continuous,   $d\rho (\lambda)= w(\lambda)d\lambda $, and     the weight   $w(\lambda)$ satisfies certain regularity conditions. The assumption $\supp \rho\subset [-1,1]$  accepted in \cite{Bern, Sz} was later partially removed in
 \cite{Gon, Nik}.
In recent years the implication $d\rho(\lambda)\to P_{n} (z)$   was successfully  developed with a help of   the Riemann-Hilbert problem  method  (see the book by P.~Deift \cite{Deift}).

We suppose that the coefficients $a_{n}, b_{n}$ are known and deduce    asymptotic properties of  the polynomials $P_n (z)$ from the behavior of $a_{n}, b_{n}$ as $n\to\infty$. Apparently, the line of research $(a_{n}, b_{n})\to P_{n} (z)$  was initiated  by P.~G.~Nevai (see his book \cite{Nev}). His approach is discussed in more details in Sect.~1.5.

    We   study the case
 \begin{equation}
\lim_{n\to\infty}\alpha_{n }  =  \lim_{n\to\infty}b_{n }  = 0 
\label{eq:comp}\end{equation}
when the perturbation $V =J-J_{0}$   is compact. Then
  the essential spectrum of the operator $J$ coincides with the interval $[-1,1]$, and its point spectrum    consists of simple eigenvalues accumulating, possibly, to the points $1$ and $-1$. We always assume that condition \e{eq:comp} is satisfied.
    Our objective is to investigate the case where $\alpha_{n}\to 0$ and $b_{n}\to 0$ as $n^{-1}$ or slower.  More precisely, we assume that
   \begin{equation}
\sum_{n=0}^\infty \big(| a_{n+1} - a_{n} |+|b_{n+1} - b_{n}|\big)<\infty .
\label{eq:LR}\end{equation}
Thus, the sequences $\{a_{n} \}$ and $\{b_{n} \}$ are of bounded variation.
In  the quantum mechanical terminology,   such perturbations $V$ of the operator $J_{0}$ are called long-range.   

 The traditional approach to 
  scattering  theory for differential operators   relies on  a study of the so called Jost solutions
$  f(x,z) $ of the Schr\"odinger equation
 \begin{equation}
 - (p (x) f' (x,z))'+q(x)f(x,z)=zf(x,z)
\label{eq:Sch}\end{equation}
  distinguished by their asymptotics as $x\to\infty$. 
  We follow the same scheme and so start with a construction of  solutions $  f(z) =\{f_{n} (z)\}_{n=-1}^\infty$   of the 
     second order Jacobi difference equation    
 \begin{equation}
 a_{n-1} f_{n-1}  (z) +b_{n} f_{n} (z) + a_{n} f_{n+1} (z)= z f_{n} (z), \q n\in{\Bbb Z}_{+}, 
\label{eq:Jy}\end{equation}
(the number $a_{-1}\neq 0$   may be chosen at our convenience; for definiteness, we put $a_{-1}=1/2$)
satisfying a certain asymptotic condition as $n\to\infty$. We   call $f(z)$ the Jost solutions and $f_{-1}  (z) $ the Jost functions. The Jost function is related to the Wronskian $\{ P(z), f(z)\}$ of the solutions $P(z)=\{ P_{n} (z)\}$ and $f(z)=\{ f_{n} (z)\}$ of equation \e{eq:Jy} by the identity
\begin{equation}
-2^{-1} f_{-1}(z)=  \{ P(z), f(z)\}=:\Omega (z).
 \label{eq:wei}\end{equation}
 We express all coefficients in asymptotic formulas for $P_{n}(z)$ in terms   of the Jost function   $f_{-1} (z)$.

  The asymptotics of the Jost solutions for the Schr\"odinger equation \e{eq:Sch} with long-range coefficients is given by the famous semiclassical  Liouville-Green Ansatz.  So, our first goal is to find its analogue for Jacobi operators. Then we develop spectral   theory of Jacobi operators $J$ with long-range coefficients essentially along the  same lines (see \cite{Y/LD}) as in the short-range case when $\alpha_{n} \to 0$ and $b_{n}\to 0$ faster than  $n^{-1}$.
We emphasize that  in the problem we consider, the semiclassical approximation   applies for large $n$ when oscillations of the coefficients $\alpha_{n}$ and $b_{n}$ are not too strong.

Of    course the asymptotic formulas we obtain are quite different for regular points $z$ of $J$, for its eigenvalues and for  $z=\lambda\in (-1,1)$ lying on its continuous spectrum. 
   Since $P_{n} (\lambda)$
 are the continuous spectrum eigenfunctions of   $J$, the last problem is natural to consider in the scattering theory framework.

\subsection{Short-range perturbations}

Let us describe the main ideas of our approach  presented in \cite{Y/LD} for
 the case of short-range    perturbations $V=  J-J_{0}$  when the condition 
 \begin{equation}
\sum_{n=0}^\infty \big(| a_{n} -1/2|+|b_{n}|\big)<\infty 
\label{eq:Tr}\end{equation} 
 is satisfied.  Then the Jost solutions of  equation \e{eq:Jy} are
  distinguished by their asymptotics
  \begin{equation}
 f_{n} (z)\sim \z (z)^n  
\label{eq:jost}\end{equation}
as $n\to\infty$. Here 
\begin{equation}
 \z (z) =z-\sqrt{z^2 -1}= (z + \sqrt{z^2 -1})^{-1}
 \label{eq:ome}\end{equation}
(we choose $\sqrt{z^2 -1}>0$ for $z>1$).   The sequence $f_{n}(z)$ rapidly tends to zero   for $z\in  {\Bbb C}\setminus [-1,1]$ and oscillates for $z=\lambda\pm i0$, $\lambda\in (-1,1)$.  At a formal level, solutions $ f_{n} (z)$ of equation \e{eq:Jy} with asymptotics \e{eq:jost} were introduced in \cite{Case}. Under assumption  \e{eq:Tr} their existence was proven in   \cite{Y/LD}. It is important that $ f_{n} (z)$ are analytic functions of $z\in{\Bbb C}\setminus [-1,1]$ and are continuous up to the cut along $(-1,1)$.

Once these results are established, spectral and scattering theories  for the operator $J$ can be developed  quite in the same way as for differential operators. Since 
     \begin{equation}
P_{n} (\lambda)=\frac{   f_{-1} (\lambda+i0)  f_{n} (\lambda-i0) -  f_{-1} (\lambda-i0)  f_{n} (\lambda+i0)}{2i\sqrt{1-\lambda^2}},  \q \lambda\in (-1,1), \q n=0,1,2, \ldots, 
\label{eq:HH4}\end{equation}
we see that $  f_{-1} (\lambda+i0) =\ov{  f_{-1} (\lambda-i0) }\neq 0$. It easily follows that  the   spectrum   of the operator $J$ is absolutely continuous on the interval $(-1,1)$ and $d\rho(\lambda)=w(\lambda) d\lambda$ with  a continuous and positive    weight  
\begin{equation}
w (\lambda)=\frac{2}{\pi}\sqrt{1-\lambda^2}|f_{-1} (\lambda+i0)|^{-2}.
 \label{eq:wei1}\end{equation}
 The   operator $J$ can also have infinite number of simple eigenvalues accumulating, possibly, to the points $1$ and $-1$; it is not excluded that these points are eigenvalues of $J$.

Relations \e{eq:jost} and \e{eq:HH4} lead  to the asymptotic formula as $n\to\infty$ for orthogonal polynomials:
 \begin{equation}
 P_{n} (\lambda)= (2/\pi)^{1/2} w(\lambda)^{-1/2}(1-\lambda^2)^{-1/4}  \sin ( (n+1)\arccos\lambda +\pi \xi(\lambda) ) + o(1)  ,  \q \lambda\in (-1,1).
\label{eq:zF}\end{equation}
The function $\xi(\lambda)$  has different names: a) In scattering theory it is known as the scattering phase or
   phase shift.   b)  It can be also identified with
  the Kre\u{\i}n spectral shift function  for the pair $J_{0}$, $J$, that is
   \[
 \xi(\lambda)= \pi^{-1}\lim_{\varepsilon\to +0}\arg\Det \big( I+ V  (J_{0}-\lambda-i \varepsilon)^{-1}\big).
\]
c) In orthogonal polynomial theory, it is constructed as the Hilbert transform of the function $\ln w(\lambda)$ and is  known as the Szeg\H{o} function. It is shown in \cite{Y/LD} that
\[
\Det \big( I+ V  (J_{0}-z)^{-1}\big)= A \z (z)  f_{-1}(z) \q \mbox{where} \q A=\prod_{k=0}^\infty(2a_{k}),
\]
so that in view \e{eq:wei1},  $\ln w(\lambda)$ and $\xi(\lambda)$ are (up to insignificant factors)    harmonic conjugate functions. Thus,  
 the spectral shift and Szeg\H{o} functions are different definitions of the same object. 
 
 Formula \e{eq:zF} is of course the same as in the classical works \cite{Bern, Sz}, but our definition of the phase shift $\xi(\lambda)$ and our assumptions on $J$ are quite different from \cite{Bern, Sz}. We emphasize that the inclusion $\supp\rho\subset [-1,1]$ is irrelevant for our approach. 
   Note that for the operator $J_{0}$, the asymptotic formula 
  \e{eq:zF} reduces to the exact expression for
the   normalized  Chebyshev polynomials   of the second kind.

 Formula \e{eq:zF} is very natural  from the scattering theory viewpoint. Indeed, the solutions $f_{n} (z)$ of the Jacobi equation  \e{eq:Jy} with asymptotics \e{eq:jost}  are discrete analogues of the Jost solutions $f(x,z)$ of the Schr\"odinger equation \e{eq:Sch}
where  $p (x)>0$,  $p (x)\to p_{0}>0$ as $x\to\infty$ and $p (x)-p_{0}$, $q(x)$ are in $L^1$. These conditions correspond to  \e{eq:Tr}.
The Jost solution  is  distinguished by its asymptotics 
  \begin{equation}
f(x,z)\sim e^{-x \sqrt{-z/ p_{0}} },   \q x\to \infty,
\label{eq:Sch1}\end{equation}
 where $\Re\sqrt{-z/ p_{0}}\geq 0$, while the regular solution  
    $  \varphi(x,z)$ of \e{eq:Sch}  is fixed by the conditions $  \varphi(0,z)=0$, $  \varphi'(0,z)=1$.
    The regular solution  plays the role of the polynomial solution $P(z)=\{ P_{n} (z)\} $ of equation \e{eq:Jy}. Since  $  \varphi(x,\lambda)$ is a linear combination of the Jost solutions $f_{n}(\lambda\pm i0)$ (cf. \e{eq:HH4}), it can be standardly (see, e.g.,   \S 4.1 and \S 4.2 of   \cite{YA})  deduced from \e{eq:Sch1} for $z=\lambda \pm i0$
that
   \begin{equation}
   \varphi(x,\lambda)=   (2/\pi)^{1/2}  w(\lambda)^{-1/2}\lambda^{-1/4}  \sin ( x\sqrt{\lambda/ p_{0}} -\pi \xi(\lambda) ) + o(1) , \q\lambda > 0, \q x\to\infty.
\label{eq:zC}\end{equation}
  Here $w(\lambda)$ is the derivative of the spectral measure of the operator $H=D p (x)D + q(x)$, and $\xi(\lambda)$  is    the spectral shift function   for the pair of the operators $H_{0}= p_{0} D^2$, $H $.  
  Obviously,   formulas \e{eq:zF} and \e{eq:zC} are quite analogous with $n$ and $x$  playing  similar roles.


 
 \subsection{Scheme of the approach}
 
 Our main goal is to study Jacobi operators $J=J_{0}  + V$ with the coefficients $\alpha_{n}$ and $b_{n}$ decaying so slowly that the condition \e{eq:Tr} is not satisfied. Instead we require a weaker assumption \e{eq:LR}.
 We follow here closely the well developed semiclassical approach in the theory of    differential equations  reviewed recently in \cite{Y-LR}. The starting point of this approach is to find a suitable modification   of the Jost solutions for the long-range case. The paper \cite{Y-LR} relied on a simplified Liouville-Green  Ansatz given by the   formula 
 \begin{equation}
f(x ,z) \sim \exp \Big(-\int_{0}^x \sqrt{\frac{q(y)-z  } {p(y)} }dy\Big), \q x\to\infty,  
\label{eq:A7}\end{equation}
for solutions of equation  \e{eq:Sch}. 

 In this paper we   combine the methods of \cite{Y/LD} where short-range perturbations of Jacobi operators were considered and of \cite{Y-LR} where differential operators with long-range coefficients were studied.   Let us  discuss the  main steps of our approach in more details:

a. 
We show in Sect.~2.2 that an analogue of \e{eq:A7} for solutions of the  difference equation \e{eq:Jy}  is given by the formula
 \begin{equation}
f_{n}(z) \sim \z \big( \frac{z-b_{0}}{2a_{0}}\big)  \z \big( \frac{z-b_{1}}{2a_{1}}\big) \cdots \z \big( \frac{z-b_{n-1}}{2a_{n-1}}  \big) =: q_{n} (z)  ,\q n \to\infty,
\label{eq:jost1}\end{equation}
where the function $\z (z)$ is defined by   \e{eq:ome}.  This means that the relative remainder
\begin{equation}
 r_{n} (z) : = q_{n} (z)^{-1} \big(a_{n-1}q_{n-1} (z) + (b_{n}-z)q_{n} (z) + a_{n}q_{n+1} (z)\big), \q n\in{\Bbb Z}_{+},
\label{eq:re2}\end{equation}
belongs to $ \ell^1 ({\Bbb Z}_{+})$. Note that, unlike \e{eq:jost},  asymptotic formula  \e{eq:jost1} takes into account decay properties of the coefficients $\alpha_{n}$ and $b_{n}$.
Since $|\z(z)|< 1$ for $z\in{\Bbb C}\setminus [-1,1]$ and $2a_{n} \to 1$, $b_{n}\to 0$, we see that the Jost solutions $f_{n}  (z)= O(\varepsilon^n)$ for some $\varepsilon= \varepsilon(z)<1$ as $n\to\infty$.  It is also easy to see that  $f_{n}  ( \lambda\pm i0)$   oscillate as $n\to\infty$ if   $\lambda\in (-1,1)$.

b.
Then we make in Sect.~2.3 a multiplicative change of variables
 \begin{equation}
   f_{n}( z )= q_{n} (z )u_{n} (z )
   \label{eq:Jost}\end{equation} 
   and reduce   the difference equation \e{eq:Jy}  to a Volterra ``integral" equation for the sequence $u_{n} (z)$
   satisfying the condition $u_{n} (z)\to 1$ as $n\to\infty$.
   
   c.
This equation is standardly solved    by iterations in Sect.~2.4 which allows us to prove  that $ u_{n} (z)$ are analytic functions of $z\in{\Bbb C}\setminus [-1,1]$ and are continuous up to the cut along $(-1,1)$. According to \e{eq:Jost} the same is true for the functions $ f_{n} (z)$ (see Theorem~\ref{EIK}).   All standard statements about spectral properties of the Jacobi operator $J$ and  asymptotic formulas for the polynomials $P_{n} (z)$ are easy consequences of   this analytic result.
   
   d. 
  If $z\in{\Bbb C}\setminus [-1,1]$ but is not an eigenvalue of the operator $J$, we prove in Theorem~\ref{GE1}  that
     \begin{equation}
\lim_{n\to\infty}  \Big(q_{n}  (z) P_{n}(z)\Big)= -\frac{\{ P(z), f(z)\}}{\sqrt{z^2-1}} .
\label{eq:HSq5}\end{equation}
To that end,  we construct an exponentially growing solution $g_{n}(z) $ of the equation \e{eq:Jy} by the formula 
 \[
 g_{n} (z)= f_{n} (z)\sum_{m=0}^n (a_{m-1} f_{m-1}(z) f_{m}(z))^{-1},\q n\in{\Bbb Z}_{+}.
\]
  Perhaps   this formula is new.
  
  e.
  Since formula \e{eq:HH4} remains true in the long-range case, we immediately obtain   the asymptotics
    \begin{equation}
 P_{n} (\lambda)=  (1-\lambda^2)^{-1/2}  |f_{-1} (\lambda + i0 ) |\sin \big( n\arccos\lambda +\Phi _{n} (\lambda) \big) + o(1)  ,  \q \lambda\in (-1,1),
\label{eq:LRz}\end{equation}
of the orthogonal polynomials.
The phase $\Phi_{n} (\lambda)$ is expressed     in terms of the phase of the Jost function $f_{-1} (\lambda + i0 )$. It depends explicitly on the coefficients $a_{n}$, $b_{n}$ and satisfies the condition
$\Phi_{n} (\lambda)=  o(n)$ as $n\to\infty$.  The precise definition  of $\Phi_{n} (\lambda)$ is given in Theorem~\ref{As}; see also formula  \e{eq:fF}.   In view of \e{eq:wei1}  the amplitude factors in \e{eq:zF} and \e{eq:LRz}  are the same.

  f.
  Our results on the Jost functions directly imply    that, under assumptions  \e{eq:comp}, \e{eq:LR}, the   spectrum of the Jacobi operator $J$ is absolutely continuous on   $(-1,1)$ and the corresponding  weight  $w (\lambda)$ is continuous and strictly positive
and can be constructed by formula  \e{eq:wei1}. This result is stated in Theorem~\ref{SF}. At  the same time, we obtain the limiting absorption principle for the operator $J$ stating that matrix elements of its resolvent $R(z)= (J-z)^{-1}$, $\Im z \neq 0$, are continuous functions of $z$ up to the cut along $(-1,1)$.


We emphasize that
 in contrast to differential operators, short-range perturbations obeying    \e{eq:Tr} are included in the class of 
long-range perturbations satisfying   \e{eq:LR}.
For example, condition \e{eq:LR} is satisfied if
\begin{equation}
a_{n} = 1/2+ \alpha n^{-r_{1}}+ \tilde{\alpha}_{n} , \q b_{n} = b n^{-r_{2}} +  \tilde{b}_{n}
\label{eq:LR1}\end{equation}
where  $\alpha, b \in {\Bbb R}$, $r_{1}, r_{2}\in (0, 1]$ and $\tilde{\alpha}_{n}, \tilde{b}_{n}\in \ell^1 ({\Bbb Z}_{+} )$. 
This reduces to \e{eq:Tr} if $\alpha = b =0$.
For Pollaczek polynomials,   relations \e{eq:LR1}  are true with $r_{1}= r_{2}=1$ and $\tilde{\alpha}_{n}=O(n^{-2})$, 
  $\tilde{b}_{n}=O(n^{-2})$. In this case the phase in formula \e{eq:zF} is essentially changed (see the book  \cite{Sz}, Section~5 in the Appendix). This resembles  the modification of the phase function for the Schr\"odinger operator with the Coulomb potential (see, e.g., formula (36,23)  in the book~\cite{LL}).

Condition \e{eq:LR} is very precise. Indeed, as shown in \cite{Nab} (see also the preceding paper \cite{Pos}),  there exist coefficients   $b_{n}$  decaying only slightly worse than $n^{-1}$ and oscillating as $n\to\infty$ such that the point spectrum of the corresponding Jacobi operator $J$ with $a_{n}= 1/2$ is dense in $[-1,1]$. In this case the limiting absorption principle for the operator $J$
does not of course hold.

 To emphasize the analogy between differential and difference operators, we often  use ``continuous" terminology (Volterra  integral equations, integration by parts, etc.) for sequences labelled by the discrete variable $n$. Below $C$, sometimes with indices,  and $c$ are different positive constants whose precise values are of no importance.

  \subsection{Related research}
  
  Now we are in a position to compare our approach with the paper \cite{Mate}  that is well known in the orthogonal polynomial literature. The results of this paper (see also the survey  \cite{Assche}) are  close to some of our results but not quite coincide with them, and the methods are essentially different.

The paper \cite{Mate} relies on specific methods of orthogonal polynomials theory.
  Apparently, the initial point of  \cite{Mate} is the relation 
  \[
 \lim_{n\to\infty} \frac{P_{n-1}(z)}{P_{n}(z)}=\z(z), \q \Im z\neq 0,
  \]
    established earlier by P.~Nevai in \cite{Nev}, Theorem~4.1.13, and
    improving one of  Poincar\'e's theorems.

  We proceed from spectral theory of Jacobi operators.
   One of important differences compared   to \cite{Mate} is the introduction and consistent use of   Jost solutions $f_{n} (z)$ of equation \e{eq:Jy} distinguished by their asymptotics \e{eq:jost1} as $n\to\infty$. Then relation \e{eq:HH4} yields asymptotic  formula \e{eq:LRz}. Our proof of asymptotic relation  \e{eq:HSq5}    is also motivated by results on differential equations.

  For the Schr\"odinger operator with short-range coefficients, the   scheme of Sect.~1.4 (except formula \e{eq:HSq5}) goes back to  the paper \cite{Jost} by R.~Jost; it is described in \S 4.1 of the book \cite{YA}. There exists another approach due to N.~Levinson (see the book \cite{CoLe}) adapted to difference equations  in the book \cite{B-L}. Probably the methods of Jost and Levinson are essentially equivalent, but,
for our purposes, it is more convenient to use the Jost method, moreover that it admits  a   direct quantum mechanical interpretation. 

We followed here the paper \cite{Y-LR} devoted to differential equations. Actually, the Ansatz \e{eq:A7}  appeared in \cite{Y-LR}  as an attempt to adjust the famous Liouville-Green Ansatz (which compared to \e{eq:A7} contains an additional pre-exponential factor) to difference equations. This allowed us in \cite{Y-LR} to get rid of conditions on the second derivatives of $p(x)$ and $q(x)$ required by the Liouville-Green Ansatz
and to develop spectral theory of the Schr\"odinger operator $H$ under the assumptions $p'\in L^1 ({\Bbb R}_{+})$, $q'\in L^1 ({\Bbb R}_{+})$.



   Let us try to advocate our approach. Its specific feature is that we treat the discrete and continuous  cases   on equal footing. Thus we can use traditional results well known for  differential operators in the discrete case. Here are some examples: 
   
   1. Both asymptotic relations \e{eq:HSq5} and \e{eq:LRz}   were obtained in \cite{Mate}. However expressions for the coefficients in the right-hand sides were not, at least in the author's opinion,  very efficient.
     It was conjectured in \cite{Mate} that the asymptotic coefficients in    \e{eq:LRz} can be obtained from that in \e{eq:HSq5} as the limit on $(-1,1)$ from complex values of $z$. This conjecture was later justified in \cite{ Va-As}. In our approach this problem does not even arise since both coefficients are expressed in terms of the Wronskian $\{P(z), f(z)\}$ of the   polynomial and Jost solutions of equation \e{eq:Jy}. 
     
     2. Introducing and consistently using the Jost solutions we avoid many more recent methods, for example, transfer matrices, Pr\"ufer variables, etc.
     
     3. The   scheme used here yields automatically an expression for the spectral measure in terms of the Jost function. This shows, in particular, that that this measure is absolutely continuous and the corresponding density is a continuous positive function. 
     
     Note that quite general conditions of the absolutely continuity of the spectral measure were given in the paper \cite{Stolz} and the book  \cite{Weidmann}, Theorem~14.25, where the subordinacy method of \cite{GD} was used.

  I thank G.~\'Swiderski who kindly informed the author about the papers \cite{Mate, Va-As}.


 \section{Volterra integral equation}
 
Here we reduce
 the Jacobi difference equation  \e{eq:Jy}  with asymptotics \e{eq:jost1}  to a Volterra equation whose solution can be constructed by iterations. Below conditions \e{eq:comp} and \e{eq:LR} are always assumed unless indicated otherwise.

\subsection{Preliminaries}

Let us consider equation \e{eq:Jy}. Note that the values of $f_{N-1}$ and $f_{N }$ for some $N\in{\Bbb Z}_{+}$ determine the whole sequence $f_{n}$ satisfying the difference equation \e{eq:Jy}.

 Let $f=\{ f_{n} \}_{n=-1}^\infty$ and $g=\{g_{n} \}_{n=-1}^\infty$ be two solutions of equation \e{eq:Jy}. A direct calculation shows that their Wronskian
  \begin{equation}
\{ f,g \}: = a_{n}  (f_{n}  g_{n+1}-f_{n+1}  g_{n})
\label{eq:Wr}\end{equation}
does not depend on $n=-1,0, 1,\ldots$. In particular, for $n=-1$ and $n=0$, we have
 \begin{equation}
\{ f,g \} = 2^{-1} (f_{-1}  g_{0}-f_{0}  g_{-1}) \q {\rm and} \q \{ f,g \} = a_{0}  (f_{0}  g_{ 1}-f_{ 1}  g_{0}).
\label{eq:Wr1}\end{equation}
Clearly, the Wronskian $\{ f,g \}=0$ if and only if the solutions $f$ and $g$ are proportional.


It is convenient to introduce a notation
\[
x_{n}'= x_{n+1}  -x_{n}
\]
for the ``derivative" of a sequence $x_{n}$. We  note the Abel summation formula (``integration by parts"):
 \begin{equation}
\sum_{n=N }^ { M} x_{n}  y_{n}' = x_{M}  y_{M+1} - x_{N -1}  y_{N}  -\sum_{n=N } ^{M} x_{n-1}'  y_{n};
\label{eq:Abel}\end{equation}
here $M\geq N\geq 0$ are arbitrary, but we have to set $x_{-1}=0$ so that $x_{-1}'=x_{0}$.

Let us  fix the branch of the analytic function $\sqrt{z^2 -1}$ of $z\in {\Bbb C}\setminus [-1,1]$ by the condition $\sqrt{z^2 -1}>0$ for $z>1$. Obviously, this function is continuous up to the cut along $[-1,1]$,    it equals $\pm i\sqrt{1-\lambda^2}$ for $z=\lambda\pm i0$, $\lambda\in (-1,1)$, and $\sqrt{z^2 -1}< 0$ for $z< -1$. Define the one-to-one, onto mapping
$\z: {\Bbb C}\setminus [-1,1] \to {\Bbb D} $ (the unit disc) by formula
 \e{eq:ome}.  
Since $2z= \z (z) +\z (z)^{-1}$,
the sequence $\{ \z (z)^n\}_{n=-1}^\infty$  satisfies the ``free" equation \e{eq:Jy} where $a_{n}=1/2, b_{n}=0$.
  For $\lambda\in [-1,1]$, it is common to set $\lambda=\cos \theta$ with $\theta\in [0,\pi ]$. Then 
$
\z (\lambda\pm i0)=e^{\mp i \theta}$.

Below the values of $|z|$ are bounded, and hence the values of $\z=\z(z)$ are separated from $0$.

   \subsection{Ansatz}

Here we show that the sequence $q_{n}  (z)$ defined in  \e{eq:jost1} is an approximate solution (Ansatz) of equation \e{eq:Jy}. 


Let $\Pi={\Bbb C}\setminus {\Bbb R}$, and let $\clos\Pi$ be the closure of $\Pi$. 
For $z\in\clos\Pi$, we set
 \begin{equation}
 z_{n}  = \frac{z-b_{n}} {2a_{n}}
\label{eq:aa}\end{equation}
and
 \begin{equation}
 \z_{n}  = \z (z_{n})=  z_{n}-\sqrt{z^2_{n} -1} = ( z_{n}+\sqrt{z^2_{n} -1})^{-1} .
\label{eq:aa1}\end{equation}
Note that
\begin{equation}
a_{n} \z_{n} + a_{n} \z_{n}^{-1}  + b_{n}-z=0.
\label{eq:re5}\end{equation}
Obviously, $z_{n} \to z$  and $\z_{n} \to \z (z)$ as $n\to\infty$ because $ 2a_{n }\to 1$ and $b_{n}\to 0$.

Now, we define  a sequence $q_{n}= q_{n} (z)$ by the relations $q_0 (z) = 1$ and  
\begin{equation}
 q_{n} (z)  = \z_{0}\z_{1} \cdots  \z_{n-1}, \q n\geq 1,
\label{eq:re1}\end{equation}
which coincides with the right-hand side of \e{eq:jost1}. 
The functions $q_{n} (z)$ are analytic in $z\in\Pi$ and   continuous up to the real axis. Of course $\z(\bar{z})=\ov {\z(z)}$ and $q_{n}(\bar{z})=\ov {q_{n}(z)}$. 
Note also that    
\begin{equation}
 |q_{n} (z)|   \leq 1,
\label{eq:Gr6a}\end{equation}
but $q_{n} (z)\neq 0$
 for all $n$ and all $z\in \clos\Pi$. We consider $ q_{n} (z) $ for $z\in \Pi$ (not for $z\in {\Bbb C}\setminus [-1,1]$)  because the values of  $\z_{n} (\lambda+ i0)$ and $\z_{n} (\lambda- i0)$ may be different 
   for $|\lambda|>1$ if $n$ is not too large.
 

 We will show
 that the sequence $ q_{n} (z)$ satisfies approximately equation \e{eq:Jy}. Let us  introduce a (relative) remainder in this equation by the formula \e{eq:re2}
(the values of $q_{-1} (z)$ and $r_{0} (z)$ are of course inessential). Since $q_{n-1} q_{n}^{-1}= \z_{n-1}^{-1}$ and $q_{n+1} q_{n}^{-1}= \z_{n} $, we can rewrite 
  \e{eq:re2} as
 \begin{equation}
   r_{n}  =  a_{n-1} \z_{n-1}^{-1} +a_{n} \z_{n} + b_{n}-z ,
 \label{eq:re2x}\end{equation}
   or using \e{eq:re5} as
\begin{equation}
 r_{n}  
 = a_{n-1} \z_{n-1}^{-1} - a_{n} \z_{n}^{-1}
 =  (a_{n-1} -a_{n})\z_{n-1}^{-1} +a_{n} (\z_{n-1}^{-1} - \z_{n}^{-1}).
\label{eq:re4}\end{equation}

\begin{lemma}\label{re}
For $z\in\clos{\Pi}$, an estimate
\begin{equation}
 | r_{n} (z) | \leq C  \frac{ |a_{n}-a_{n-1}| + |b_{n}-b_{n-1}|}
 {\sqrt{z_{n-1}^2-1}+\sqrt{z_{n}^2-1}} 
\label{eq:re3}\end{equation}
holds.
In particular, $\{ r_{n} (z) \}_{n=0}^\infty \in \ell^1 ({\Bbb Z}_{+})$ if $z\neq\pm 1$.
\end{lemma}
    
\begin{pf}
Let us proceed from representation \e{eq:re4}.
 By definition  \e{eq:aa1}, we have 
 \begin{equation}
\z_{n-1}^{-1} - \z_{n}^{-1}= (z_{n-1} - z_{n})\Big(1+  \frac{z_{n-1}+ z_{n}}{\sqrt{z_{n-1}^2-1}+\sqrt{z_{n}^2-1}}\Big),
\label{eq:aa1-}\end{equation}
where according to  \e{eq:aa}
\[
 z_{n-1} - z_{n}= \frac{(z-b_{n} )(a_{n}-a_{n-1}) + a_{n }(b_{n}-b_{n-1}) }{2 a_{n-1}a_{n}} 
\]
so that 
 \begin{equation}
|z_{n} - z_{n-1}| \leq C  (|a_{n} - a_{n-1}| +|b_{n} - b_{n-1}| ) . 
\label{eq:aax}\end{equation}
Thus equality \e{eq:aa1-} yields estimate \e{eq:re3}.
\end{pf}

\subsection{Multiplicative substitution}

Let the sequence
$ q_{n} (z)$ be given by formulas   \e{eq:aa}, \e{eq:aa1} and \e{eq:re1}.
We are looking for solutions 
  $f_{n} (z)$  of the difference equation \e{eq:Jy} satisfying the condition 
   \begin{equation}
f_{n} (z) = q_{n} (z)(1+ o(1)), \q n\to\infty.
\label{eq:Gs1}\end{equation}
The uniqueness of such solutions is almost obvious.

\begin{lemma}\label{uniq}
Equation \e{eq:Jy} may have only one solution   $f_{n} (z)$ satisfying condition \e{eq:Gs1}.
 \end{lemma}

\begin{pf}
Let $\tilde{f}_{n} (z)$ be another solution of \e{eq:Jy} satisfying  \e{eq:Gs1}. Then the Wronskian \e{eq:Wr} of these solutions
equals
\[
\{f,\tilde{f}\}=a_{n} q_{n} q_{n+1} o(1)  .
\]
This expression tends to zero because the sequences $a_{n} $ and $q_{n} $ are bounded according to  \e{eq:comp} and \e{eq:Gr6a}. Thus $\tilde{f} =Cf$ where $C=1$ by virtue again of condition \e{eq:Gs1}.
\end{pf}

 
For construction of $f_{n}  (z)$, we will reformulate the problem    introducing a  sequence
\begin{equation}
 u_{n} (z)=  q_{n} (z)^{-1}  f_{n} (z), \q n\in {\Bbb Z}_{+}.
\label{eq:Gs4}\end{equation}
 Then \e{eq:Gs1} is equivalent to the condition  
\begin{equation}
\lim_{n\to\infty} u_{n} (z)=   1.
\label{eq:A12a}\end{equation}

Let us derive a difference equation for $ u_{n} (z)$.

\begin{lemma}\label{subst}
Let $z\in\clos\Pi$ and let $r_{n} (z) $ be given by formula \e{eq:re2}. Then
 equation  \e{eq:Jy} for a sequence $ f_{n} (z)$ is equivalent to the equation
\begin{equation}
 a_{n} \z_{n} ( u_{n+1} (z)- u_{n} (z)) - a_{n-1} \z_{n-1}^{-1} ( u_{n} (z)- u_{n-1} (z))=- r_{n} (z) u_{n} (z), \q n\in {\Bbb Z}_{+},
\label{eq:A12}\end{equation}
for  sequence  \e{eq:Gs4}. 
 \end{lemma}

\begin{pf}
Substituting  expression $f_{n}  = q_{n}  u_{n} $ into \e{eq:Jy}, we see that
\begin{align*}
q_{n}^{-1}\Big(  a_{n-1} f_{n-1} + ( b_{n} -z)f_{n} + a_{n} f_{n+1}\Big)
 =  & a_{n-1} \z_{n-1}^{-1} u_{n-1} + ( b_{n} -z)u_{n} + a_{n} \z_{n} u_{n+1}  
\\
      & a_{n-1} \z_{n-1}^{-1} (u_{n-1}-u_{n}) + r_{n} u_{n} + a_{n} \z_{n} (u_{n+1} -u_{n}) .
\end{align*}
  by virtue of \e{eq:re2x}. Thus equations \e{eq:Jy}  and \e{eq:A12} coincide.
\end{pf}

According to Lemma~\ref{re} the sequence $r_{n} (z) $ of the coefficients of $u_{n} (z)$ in the right-hand side of   \e{eq:A12} belongs to $  \ell^1 ({\Bbb Z}_{+})$. This allows us  to reduce the  difference equation \e{eq:A12} with condition \e{eq:A12a}   to a  ``Volterra integral" equation 
\begin{equation}
 u_{n}(  z )= 1 -  \sum_{m=n+1}^\infty G_{n,m} (z)     r_{m} ( z )u_{m} (  z )   
\label{eq:A17}\end{equation}
 with kernel 
\begin{equation}
G_{n,m} (z) =  q_{m}(z)^2 \sum_{p=n}^{m-1} ( a_{p} \z_{p })^{-1}  q_{p }(z)^{-2},\q n,m\in{\Bbb Z}_{+} , \q m\geq n+1.
\label{eq:Gg1}\end{equation}
Note that $G_{n,m} (\bar{z})=\ov{G_{n,m} (z)}$.
The functions $G_{n,m} (z) $ are analytic in $z\in\Pi$ and are continuous up to the real axis.

\subsection{Integral equation}

Our plan is now the following. We first prove that, for all $z\in \clos \Pi\setminus \{-1,1\}$, a solution $u_{n} (z) $ of the integral equation \e{eq:A17} exists for sufficiently large $n$ and tends to $1$  as $n\to\infty$. Then we 
show that, for such $n$, the sequence $u_{n} (z) $ satisfies also the difference equation \e{eq:A12}.
 To all $n\geq -1 $, the sequence $u_{n} (z) $  is extended as a solution of equation \e{eq:A12}.


The following assertion plays the crucial role in our analysis of equation \e{eq:A17}, in particular, for $z$ lying on the cut along $[-1,1]$. It shows that  the sequence \e{eq:Gg1}  is bounded   uniformly in $n$ and $m$  provided the points $ \pm 1$ are excluded.

 \begin{lemma}\label{eik1}
There exist constants $C(z)$ and $N(z)$ such that  an estimate 
\begin{equation}
| G_{n,m} (z) |\leq  C (z)<\infty, \q m-1\geq n \geq N(z),
\label{eq:A19}\end{equation}
 is true for all $z\in \clos\Pi\setminus \{-1,1\}$. The constants $C(z)$ and $N(z)$ are common for $z$ in  compact subsets of $\clos\Pi\setminus \{-1,1\}$, that is, for all
$ z\in \clos\Pi$ such that $|z^2-1|\geq \epsilon$ and $|z|\leq R$ where   $\epsilon>0$ and $R<\infty$ are some fixed numbers. 
 \end{lemma}

 \begin{pf}
   According to definition \e{eq:re1} we have
 \[
   (q_{p}^{-2})' = ( \z_{p}^{-2} -1)q_{p}^{-2} .
\]
Set 
\[
\eta_{p}=  a_{p}^{-1}  \z_{p} (1-\z_{p}^2)^{-1}. 
\]
According to \e{eq:comp}
    $\z_p\to \z$ where $\z^2\neq 1$ as $p\to\infty$, so that   $|\eta_{p}|\leq C<\infty $ for $ p\geq N $   if $N$ is sufficiently large. Therefore
 it follows (cf. \e{eq:aa1-} and \e{eq:aax}) from condition \e{eq:LR}   that  
  \begin{equation}
 \q\sum_{p=N}^\infty |\eta_{p}'|   \leq C   \sum_{p=N}^\infty (|a_p' | +|  \z_p'| )
  \leq C_{1}   \sum_{p=N}^\infty (|a_p' | +|  b_p'| ) <\infty.
 \label{eq:AXY}\end{equation}
   Integrating by parts, that is, using identity \e{eq:Abel}, we find that
 \[
 \sum_{p=n}^{m-1} (a_{p}\z_{p})^{-1} q_{p}^{-2}  =
    \sum_{p=n}^{m-1}  \eta_{p} ( q_{p}^{-2})'
 = \eta_{m}  q_{m+1}^{-2} -\eta_{n-1} q_{n}^{-2} 
 -\sum_{p=n}^{m-1}\eta_{p-1} ' q_{p}^{-2}.
\]
 Substituting this expression into \e{eq:Gg1} and using the estimates $|q_{m} q_{p}^{-1}|\leq 1$ for $m\geq p$ and
 \e{eq:AXY}, we see that 
  \[
| G_{n,m} |\leq C \big( 1+
 \sum_{p=n}^{m-1}|\eta_{p-1} '|   |q_{m}^2 q_{p}^{-2}|  \big) \leq C_{1}<\infty.
 \]
 This proves \e{eq:A19}.
  \end{pf}

  Lemmas~\ref{re} and \ref{eik1} allow us to solve the Volterra equation \e{eq:A17} by iterations.

  \begin{lemma}\label{Jos2}
   Let  $z\in \clos \Pi\setminus \{-1,1\}$. Set $u^{(0)}_n (z)=1$   and 
  \begin{equation}
 u^{(k+1)}_{n}(z)=- \sum_{m=n +1}^\infty G_{n,m} (z)  r_{m} (z)u^{(k )}_m(z),\q k\geq 0,
\label{eq:V5}\end{equation}
for all $n\in {\Bbb Z}_{+}$.
Then estimates
 \begin{equation}
| u^{(k )}_{n}(z)|\leq \frac{C(z)^k}{k!}\Big(\sum_{m=n+1 }^\infty |r_{m}(z)|\Big)^k,\q \forall k\in{\Bbb Z}_{+},
\label{eq:V6}\end{equation}
are true for all sufficiently large $n $ with the same constant $C(z) $ as in Lemma~\ref{eik1}. 
\end{lemma}

  \begin{pf}
  Suppose that \e{eq:V6} is satisfied for some $k\in{\Bbb Z}_{+}$. We have to check 
 the same estimate (with $k$ replaced by $k+1$ in the right-hand side)  for $ u^{(k+1)}_{n}$.  
 Set
 \[
 {\cal R}_{m}= \sum_{p=m +1}^\infty  | r_p | .
 \]
  According to definition \e{eq:V5}, it  follows from  \e{eq:A19} and  \e{eq:V6} that
   \begin{equation}
| u^{(k +1)}_{n} |\leq \frac{C^{k+1}}{k!}  \sum_{m=n +1}^\infty  | r_m | {\cal R}_{m}^k.
\label{eq:V7}\end{equation}
Observe that
 \[
{\cal R}_{m}^{k+1}+ (k+1)   | r_{m}|  {\cal R}_{m}^k  
\leq
{\cal R}_{m-1}^{k+1},
\]
and hence, for all $N\in{\Bbb Z}_{+}$,
   \[
 (k+1)  \sum_{m=n +1}^N  | r_m | {\cal R}_{m}^k 
 \leq 
 \sum_{m=n +1}^N  ( {\cal R}_{m-1}^{k+1}-{\cal R}_{m}^{k+1})
= {\cal R}_{n}^{k+1}-{\cal R}_{N}^{k+1}\leq  {\cal R}_{n}^{k+1}.
 \]
Substituting this bound into   \e{eq:V7}, we obtain estimate \e{eq:V6} for $u^{(k +1)}_{n}$.
    \end{pf}
  
Now we can conclude our study of the ``integral"  equation \e{eq:A17}.

  \begin{theorem}\label{eik2}
 Let assumptions \e{eq:comp} and \e{eq:LR} be satisfied. 
  For  $z\in \clos \Pi\setminus\{-1,1\}$,
equation  \e{eq:A17} has a $($unique$)$ bounded solution $\{u_{n}( z )\}_{n=0}^\infty$. This sequence obeys an estimate
 \begin{equation}
| u_{n}( z )-1|\leq C \varepsilon_{n}
\label{eq:A20}\end{equation} 
where the constant $C$ is common for $z$ in  compact subsets of $\clos\Pi\setminus\{-1,1\}$ and
  \begin{equation}
\varepsilon_{n} : = \sum_{m=n}^\infty (|\alpha_{m}'| + |b_{m}'| ).
\label{eq:eps}\end{equation} 
For all $n\in {\Bbb Z}_{+}$, the functions $u_{n}( z )$ are analytic in $z\in  \Pi$  and are continuous up to the cut along $\Bbb R$
with possible exception of the points $z=-1$ and $z=1$.
  \end{theorem}
   
\begin{pf}
 Set
     \begin{equation}
   u_{n} =\sum_{k=0}^\infty u^{(k)}_{n} 
\label{eq:A2M}\end{equation}
where $u^{(k)}_{n}$ are defined by recurrence relations \e{eq:V5}.
Estimate \e{eq:V6} shows that this series is absolutely convergent. Using the Fubini theorem to interchange the order of summations in $m$ and $p$, we see that
\[
   \sum_{m=n+1}^\infty G_{n,m}       r_{m}  u_{m}    =   \sum_{k=0}^\infty\sum_{m=n+1}^\infty G_{n,m}     r_{m}  u_{m}^{(k)} =- \sum_{k=0}^\infty  u_{n}^{(k+1)}=1- \sum_{k=0}^\infty  u_{n}^{(k)}=1- u_{n}.
\]
This is equation  \e{eq:A17} for sequence \e{eq:A2M}. It also follows from \e{eq:V6}  that
\[
 |   u_{n}(z)-1 | \leq C\sum_{m=n +1}^\infty | r_{m} (z)|
\]
which in view of  \e{eq:re3} implies  \e{eq:A20}. Since every function $u^{(k)} _{n}(z)$ is analytic in $z\in\Pi$ and is continuous up to the cut $\Bbb R$ (away from the points $\pm 1$), estimate \e{eq:V6} guarantees that function \e{eq:A2M}  possesses the same properties.
\end{pf} 
 
Let us  come back to the difference equation \e{eq:A12}.
 
  \begin{lemma}\label{int1}
  For   $z\in \clos \Pi\setminus\{-1,1\}$, a solution $u_{n} (z )$ of   integral  equation  \e{eq:A17} satisfies also difference equation \e{eq:A12}.
 \end{lemma}
 
  \begin{pf}
Below we consistently take into account that $q_{n+1}=\z_{n}q_{n}$.
  It follows from \e{eq:A17}    that
  \begin{equation}
 u_{n+1} - u_{n}= -  \sum_{m=n+2}^\infty (G_{n+1, m}-G_{n, m})r_{m}u_{m}+
 G_{n, n+1} r_{n+1}u_{n+1}.
  \label{eq:A17a}\end{equation}
  Since according to  \e{eq:Gg1}  
\[
G_{n+1, m}-G_{n, m}=-  (a_{n}\z_{n})^{-1}q_{n}^{-2}  q_{m}^2\q \mbox{and}\q
G_{n, n+1}= \z_{n} a_{n}^{-1} ,
\]
equality \e{eq:A17a} can be rewritten as
\begin{equation}
a_{n}   ( u_{n+1} - u_{n}) = \z_{n}^{-1}q_{n}^{-2} \sum_{m=n+1}^\infty     q_{m}^2 r_{m} u_{m} .
  \label{eq:A17b}\end{equation}
Putting together this equality with the same equality where $n$ is replaced by $n-1$, we arrive to  equation \e{eq:A12}.
   \end{pf}
   
     \begin{remark}\label{u'}
     Let $z\in\Pi$ or $z=\lambda\pm i0$ where $|\lambda|>1$.
      Since $|\z| <1$ and $\z_{n} \to \z$ as $n\to\infty$,   we see that $|\z_{m}|\leq \varepsilon$ for some $\varepsilon< 1$ and  sufficiently large $m $ whence
\begin{equation}
 \big|  q_{n}/q_{m}\big| =| \z_{m}\cdots \z_{n-1}  |\leq \varepsilon^{n-m}, \q n\geq m.
 \label{eq:GXx}\end{equation}
  Note also that $\{r_m u_m\}\in \ell^1 ({\Bbb Z}_{+})$. Thus, it follows from representation \e{eq:A17b}   that $\{u'_{n}\}\in \ell^1 ({\Bbb Z}_{+})$.
 \end{remark}

    \section{Modified Jost solutions}
    
   \subsection{Construction}
   
      Let us put together the results obtained in the previous section. According to  Theorem~\ref{eik2} for every $z\in\clos\Pi\setminus\{-1,1\}$ there exists a solution  $u_{n}( z )$, $n\in {\Bbb Z}_{+}$, of the integral equation \e{eq:A17}.  By  Lemma~\ref{int1} it satisfies also the difference equation \e{eq:A12}. Then Lemma~\ref{subst} implies that the function
      $ f_{n}( z ):= q_{n} (z )u_{n} (z )$
   satisfies equation \e{eq:Jy}. Estimate  \e{eq:A20} for  $u_{n}( z )$ is obviously equivalent to the asymptotics
   \begin{equation}
f_{n}(  z )=q_{n}( z )\big(1 + O( \varepsilon_{n})\big), \q n\to\infty,
\label{eq:A22}\end{equation} 
for  $f_{n}( z )$.
    
    Thus we arrive at the following result.
    
    \begin{theorem}\label{EIK}
   Let assumptions \e{eq:comp} and \e{eq:LR} be satisfied, and let  $z\in \clos\Pi  \setminus \{-1,1\}$.
Denote by $u_{n}( z )$   the sequence constructed in Theorem~\ref{eik2}. Then  the sequence $f_{n}( z )$ defined by equality \e{eq:Jost} satisfies equation \e{eq:Jy},  and it has  asymptotics
\e{eq:A22}. For all $n\in {\Bbb Z}_{+}$, the functions $f_{n}( z )$ are analytic in $z\in  \Pi$  and are continuous up to the cut along $\Bbb R$
with possible exception of the points $z=-1$ and $z=1$. Asymptotics \e{eq:A22} is  uniform in $z$ from compact subsets    of the set $  \clos\Pi\setminus\{-1,1\} $.
 \end{theorem}

   Recall that the polynomials   $P_{n}(z)$  are solutions of equation  \e{eq:Jy} satisfying the conditions  $P_{-1}(z)=0$,  $P_{0}(z)=1$.
   Put $ P(z)=\{ P_{n}(z)\}_{n=-1}^\infty$, $ f(z)=\{ f_{n}(z)\}_{n=-1}^\infty$. Then the first formula \e{eq:Wr1} yields relation \e{eq:wei} for the Wronskian $\Omega(z)$ of the solutions $ P(z) $ and  $ f(z)$.
By analogy with the continuous case,   the sequence $\{f_{n}( z )\}_{n=-1}^\infty$   will be called
  the  (modified) Jost solution of equation \e{eq:Jy} and 
 $f_{-1}(z)$ will be called
     the (modified) Jost function.  For the operator $H_{0}$,     the Jost solution is $\{\z(z)^n \}_{n=-1}^\infty$ and  the corresponding Wronskian
          \begin{equation}
     \Omega_{0} (z) = -(2\z (z) )^{-1}.
     \label{eq:WRfr}\end{equation}
     
     The following result is a direct consequence of Theorem~\ref{EIK}.

\begin{corollary}\label{JOST}
The Wronskian $\Omega(z)$   depends analytically on $z\in\Pi$, and it is continuous in $z$ up to the cut along $\Bbb R$ except, possibly, the points $\pm 1$. 
\end{corollary}

 Note  that
  \[
f_{n}(  \bar{z} )=\ov{f_{n}( z )}
 \]
 and, in particular, 
  \begin{equation}
 f_{n}( \lambda- i0 )=\ov{f_{n}(  \lambda+ i0 )}, \q \lambda\in {\Bbb R}\setminus\{-1,1\}.
 \label{eq:Acx}\end{equation}
  
   \begin{remark}\label{EIKrem}
  \begin{enumerate}[{\rm(i)}]
  \item
  In contrast to the short-range case, the   functions $f_{n}  (z)$ are not analytic on the whole  set ${\Bbb C}\setminus [-1,1]$ because, in general, $f_{n}( \lambda- i0 )\neq f_{n}(  \lambda+ i0 )$ even for $|\lambda|>1$. This circumstance is, however, inessential -- see Remark~\ref{AR}  below.

\item
The  definitions of the function $q_{n}(z)$ and of the Jost solution $f_{n}(z)$       are not   intrinsic. One can multiply them by a factor depending on $z$ but not on $n$. For example,  we can set
  \[
  \tilde{f}_{n} (z) = q_{N_{0}}(z)^{-1}f_{n}(z)
  \]
  where $N_{0}$ is some fixed number. Then $  \tilde{f}_{n} (z) $ satisfies equation \e{eq:Jy} and 
  $  \tilde{f}_{n} (z) = \tilde{q}_n (z) (1+ O(\varepsilon_{n}))$ where $\tilde{q}_n(z)= q_{N_{0}}(z)^{-1} q_{n} (z)$ as $n\to\infty$.
  \end{enumerate}
 \end{remark}

     Let us find asymptotics of the sequence \e{eq:re1} as $n\to\infty$.
     
      \begin{lemma}\label{asq}
   Let assumption \e{eq:comp}  be satisfied,  let $z\in\clos\Pi\setminus\{-1,1\}$,  and let $\z=\z (z)$ be given be equality
   \e{eq:ome}. Then 
     \begin{equation}
q_{n} (z)= e^{n(\ln\z + o (1))}, \q n\to\infty.
 \label{eq:asq}\end{equation} 
 \end{lemma}
  
   \begin{pf}
   Let the sequences $z_{n}$ and $\z_{n}$ be defined by formulas  \e{eq:aa} and \e{eq:aa1}.
   It follows from \e{eq:comp}  that $z_{n}=z+ o(1)$ and hence
   $
\z_{n}  = \z (1 + \epsilon_{n})  
$
 where $\epsilon_{n}\to 0$ as $ n\to\infty$. For the sequence \e{eq:re1}, this yields
    \[
\ln q_{n}  =n \ln \z+ \sum_{m=0}^{n-1} \ln (1 + \epsilon_m )  = n \ln \z+ o(n)
 \]
 which is equivalent to \e{eq:asq}.
        \end{pf}
        
        Note that $|\z (z)| <1$ for $z\in \Pi$ and  for $z=\lambda\pm i0$ if $|\lambda| >1$.   Therefore for such $z$ according to
  \e{eq:A22} and \e{eq:asq}, $ f_{n}(z)= O(\d^n)$ with some $\d<1$ as $n\to \infty$ whence $f(z)\in \ell^2 ({\Bbb Z}_{+})$.

  

 Observe that    $\Omega(z)\neq 0$ for $z\in\Pi$. Indeed, if $\Omega(z)=0$, then, by definition \e{eq:wei}, $P_{n} (z)= c f_{n} (z)$  for some  $c\in{\Bbb C}$. Since $f(z)\in \ell^2 ({\Bbb Z}_{+})$, it follows that the complex number $z$ is 
  an eigenvalue of the self-adjoint operator $J$ which is impossible. This argument also shows that a number $\lambda\in (-\infty,-1)\cup (1,\infty)$ is an eigenvalue of  $J$ if and only if $\Omega (\lambda\pm i0)=0$. By \e{eq:Acx},    these equalities for the signs $``+"$ and $``-"$ are equivalent to each other.

\subsection{On the cut}

Suppose now that $z=\lambda\pm i0$ where $\lambda\in (-1,1)$ so that $\lambda=\cos\theta$, $\theta\in (0,\pi)$.
Set
\begin{equation}
\lambda_{n} : =   \frac{\lambda-b_{n}}{2a_{n}}.  \label{eq:Joz}\end{equation}
For $|\lambda_{n}|< 1$, we have
\[
\z (\lambda_{n} \pm i0)= \lambda_{n} \mp i \sqrt{1-\lambda_{n}^2}= e^{\mp i\theta_{n} (\lambda)}
\]
where  
 \begin{equation}
   \theta_n (\lambda)= \arccos  \lambda_{n}\in (0,\pi).
\label{eq:Joa}\end{equation}
Set $\theta_{n}(\lambda ) = 0$ for
  $\lambda_{n}\geq 1$ when  $\z(\lambda_{n})>0$ and, similarly,   $\theta_{n}(\lambda)=\pi$
for  $\lambda_{n}\leq -1$, when $\z(\lambda_{n})<0$. 
Then formula \e{eq:re1} reads as 
\[
q_{n} (\lambda \pm i0)= e^{\mp i \varphi_{n}  (\lambda)}  \prod_{m=0}^{n-1}|\z(\lambda_{m}\pm i0)|
\]
where
   \begin{equation}
  \varphi_{n}(\lambda )= \sum_{m=0}  ^{n-1} \theta_m (\lambda).
\label{eq:Jo1}\end{equation}
Since $|\z(\lambda_{m}\pm i0)|=1$ for sufficiently large $m$,  the product
   \[
k(\lambda):= \prod_{m=0}^\infty|\z(\lambda_{m}\pm i0)|
\]
consists of a finite number of terms. Obviously, $k(\lambda)$ is a continuous function of $\lambda\in (-1,1)$ and $k(\lambda)\neq 0$.

As mentioned in Remark~\ref{EIKrem} (ii), the definitions of the Jost solution $f_{n} (z)$ and hence of the  Wronskian 
$\Omega(z)$ are not unique. For $\lambda\in (-1,1)$, it is convenient to normalize them dividing  by the factor $k(\lambda)$: 
  \begin{equation}
\mathbf{f}_{n} (\lambda\pm i0)= k(\lambda)^{-1}f_{n} (\lambda\pm i0) , \q \boldsymbol{\Omega} (\lambda\pm i0)=\{P(\lambda), \mathbf{f} (\lambda\pm i0)\}=k(\lambda)^{-1}\Omega(\lambda\pm i0).
\label{eq:nJ}\end{equation}
This simplifies formulas  below.

   The following   direct consequence of  Theorem~\ref{EIK} is stated in terms of the normalized Jost solutions.
 
 \begin{theorem}\label{As}
 Let assumptions~\e{eq:comp} and \e{eq:LR} be satisfied.  For $ \lambda\in (-1,1)$, define the phases $\varphi_{n}(\lambda)$  
   by formulas \e{eq:Joz} -- \e{eq:Jo1}.
 Then
   \begin{equation}
  \mathbf{f}_{n}(\lambda\pm i0)=   e^{\mp i \varphi_{n}(\lambda)}  (1+O(\varepsilon_{n})), \q n\to\infty,
   \label{eq:Jost1}\end{equation}
 where   $ \varepsilon_{n }$ is given by \e{eq:eps}.
 \end{theorem}
 
 By definitions  \e{eq:Joz},  \e{eq:Joa}, we have 
\begin{equation}
 \lambda_{n}= \lambda+ o(1)
 \q\mbox{and}  
 \q\theta_{n}= \theta+ o (1) 
  \label{eq:JostX}\end{equation}
  where $\theta=\arccos\lambda$.
  It follows that
 \[
 \varphi_{n}(\lambda)=n\arccos\lambda + \Phi_{n}(\lambda)
 \]
 where
 \begin{equation}
 \Phi_{n}(\lambda)=\sum_{m=0}^{n-1} \big( \theta_{m}  (\lambda)-\theta\big)= o(n).
  \label{eq:fF}\end{equation}
In particular, we see that
 asymptotics \e{eq:Jost1} of  $   \mathbf{f}_{n}(\lambda\pm i0)$ as $n\to\infty$  is oscillating.  
 

\section{Orthogonal polynomials}


\subsection{Uniform estimates} 

Let us start with an  estimate on orthogonal polynomials $P_{n} (z)$ for large $n$. This estimate will be uniform in $z$ from compact subsets of ${\Bbb C}\setminus\{-1, 1\}$. We recall that $P_{n} (z)$ satisfy the difference equation \e{eq:pol}
and the  conditions $P_{-1 } (z) =0$, $P_0 (z) =1$.  As before, we define the sequence $q_{n}(z)$   by formula \e{eq:re1} but instead of \e{eq:Gs4} we make a substitution
\begin{equation}
 {\sf u}_{n} (z)=  q_{n}(z)  P_{n} (z).
\label{eq:PP}\end{equation}
Observe that $q_{n}(z)^{-1}$ is an approximate solution of equation \e{eq:Jy}, that is 
\[
{\sf r}_{n}(z) : = q_{n} (z)  \big(a_{n-1}q_{n-1} (z)^{-1} + (b_{n}-z)q_{n} (z)^{-1} + a_{n}q_{n+1} (z)^{-1}\big), \q n\in{\Bbb Z}_{+},
\]
 where   the   remainder ${\sf r}_{n} (z) $ is given  (cf. \e{eq:re2}, \e{eq:re2x}, \e{eq:re4})  by the formula  
   \[
      {\sf r}_{n}  =  a_{n-1} \z_{n-1}  +a_{n} \z_{n}^{-1} + b_{n}-z = a_{n-1} \z_{n-1}  -a_{n} \z_{n} .
\]
 Similarly to Lemma~\ref{re}, it is easy to show that $\{   {\sf r}_{n}\}\in \ell^1 ({\Bbb Z}_{+})$ and,
similarly to Lemma~\ref{subst}, it is easy to obtain an equation
 \begin{equation}
 a_{n} \z_{n}^{-1} (  {\sf u}_{n+1} (z)-  {\sf u}_{n} (z)) - a_{n-1} \z_{n-1}  (  {\sf u}_{n} (z)-  {\sf u}_{n-1} (z))=- {\sf r}_{n} (z)  {\sf u}_{n} (z), \q n\in {\Bbb Z}_{+},
\label{eq:reg1}\end{equation}
for  sequence  \e{eq:PP}. This difference equation can be reduced 
to a Volterra integral equation.

\begin{lemma}\label{int}
The sequence \e{eq:PP} obeys   an equation
   \begin{equation}
 \mathsf u_{n+1}  ( z )=  \mathsf u_{n}^{(0)}  ( z ) - \sum_{m=1}^{n} K_{n,m} (z)  {\sf r}_{m} ( z ) \mathsf u_m  ( z ),\q n\geq 1,
\label{eq:A17m}\end{equation}
where $\mathsf u_1 (z) = q_{1}  P_{1} (z)=2z_{0}\z_{0}$,
\[
\mathsf u_{n}^{(0)}  ( z )=\mathsf u_1 (z) + a_{0}\z_{0}\sum_{p=1}^{n} a_{p} ^{-1}\z_{p} q_{p}(z)^2
\]
and
\[
K_{n,m} (z) =q_{m}(z)^{-2}\sum_{p=m}^{n} a_{p} ^{-1}\z_{p} q_{p}(z)^2.
\]
 \end{lemma}

\begin{pf}
 Set
\begin{equation}
v_{n} = a_{n} (\mathsf {u }_{n+1}  -\mathsf {u }_{n} )
\q {\rm and} \q \varrho_{n}  = {\sf r}_{n}  {\sf u}_{n}  .
\label{eq:A13}\end{equation}
Then  the second order difference equation \e{eq:reg1} for $\mathsf u_{n}$ yields a first order difference equation 
\begin{equation}
   \z_{n}^{-1} v _{n} -   \z_{n-1} v _{n-1}  = - \varrho_{n}  
\label{eq:A13X}\end{equation}
  for the sequence $v_{n} $. Note that  $\z_{n} q_{n-1}^2$ is the solution of the corresponding homogeneous equation because  $q_{n}=\z_{n-1} q_{n-1}$,
  by \e{eq:re1}.  Thus, setting $v_{n}=\z_{n}  q_{n}^{2}w_{n}$, we rewrite \e{eq:A13X} as an equation
\[
w_{n} - w_{n-1}= -   q_{n}^{-2} \varrho_{n}
\]
for $w_{n}$ whence
\[
  w_{n}=w_{0}-    \sum_{m=1}^n    q _m^{-2}\varrho_m.
\]
It follows that the solution of equation \e{eq:A13X} is given by the formula
\[
  v_{n}= v_{0}  q^2_{n} \z_{n}  \z_{0}^{-1} -  q^2_{n} \z_{n} \sum_{m=1}^n   q _m^{-2}\varrho_m,\q n\geq 1.
\]
Using now \e{eq:A13}, we find that
\[
{\sf u}_{n+1}-  {\sf u}_{1}=\sum_{p=1}^n a_{p}^{-1} v_{p}=  v_{0}  \z_{0}^{-1}\sum_{p=1}^n a_{p}^{-1} q^2_p \z_p- 
\sum_{p=1}^n a_{p}^{-1}    q^2_{p} \z_{p} \sum_{m=1}^p   q _m^{-2}{\sf r}_{m}   {\sf u}_{m}.
\]
  Interchanging the       summations over $p$ and $m$ here and observing that $ v_{0}  \z_{0}^{-1}=a_{0}\z_{0}$, we obtain equation \e{eq:A17m}.  
\end{pf}

Similarly to Lemma~\ref{eik1},   integration by  parts  shows that the functions $ \mathsf u_{n}^{(0)}  ( z )$ and $K_{n,m} (z)$ are uniformly bounded. Therefore solving  equation \e{eq:A17m} by iterations (cf. Lemma~\ref{Jos2}), we see that its solution $ \mathsf u_{n}  ( z )$ is also bounded. Coming back to relation \e{eq:PP}, we can state the following result.

 \begin{theorem}\label{REGB}
   Let assumptions \e{eq:comp} and \e{eq:LR} be satisfied, and let $z\in{\Bbb C}\setminus\{-1, 1\}$.
   Define the sequence $q_{n}(z)$   by formula \e{eq:re1}. Then the orthogonal polynomials $P_{n} (z)$
  obey an estimate
     \[
 | P_{n}(z)|\leq C |  q_{n}(  z) |^{-1}, \q\forall n \in{\Bbb Z}_{+},
\]
where $C$ does not depend on $z$ in  compact subsets of ${\Bbb C}\setminus\{-1, 1\}$.
 \end{theorem}
 

 Observe that  \e{eq:A17m} is different from the standard (see, e.g., equation (3.2) in \cite{Y/LD}) equation for orthogonal polynomials relying on the perturbation theory, that is, on a comparison of $P_{n} (z)$ with the Chebyshev polynomials of the second kind. We emphasize that Theorem~\ref{REGB} is  very close to Lemma~4 in \cite{Va-As}, and we state it mainly for the   completeness of our presentation.


\subsection{Asymptotics in the complex plane}

 Here we   find    asymptotics  of the polynomials $P_{n} (z)$  for $z\not\in [-1,1]$. We follow    the scheme exposed in \cite{Y/LD} for Jacobi operators in the short-range case and  in \cite{Y-LR} for differential  operators with long-range coefficients. In this subsection we use Theorem~\ref{EIK} for a fixed $z\not\in [-1,1]$ only.

  Let the sequence  $q_{n}  (z)$ be defined by formula \e{eq:re1}, and let $f_{n}  (z)$ be the Jost solution of equation  \e{eq:Jy}.     If $z=\lambda\in {\Bbb R}  \setminus [-1,1]$, we can choose  the  sequences $q_{n}  (\lambda\pm i0)$ and $f_{n}  (\lambda\pm i0)$ for any of the signs $``\pm"$.
   We start by introducing a solution $g_{n} (z)$ of equation  \e{eq:Jy} exponentially growing as $n\to\infty$. Perhaps this construction is of interest in its own sake. 
  In view of asymptotics  \e{eq:A22}, we can choose $n_{0}= n_{0} (z)$ such that $f_{n} (z)\neq 0$ for all $n\geq n_{0}-1$. Note that, for $\Im z\neq 0$,    one can set $n_{0}= 0$ because the equality $f_{n_{0}-1} (z)= 0$  implies that the Jacobi operator $J^{(n_0)}$ with the matrix elements $a^{(n_0)}_{n}=  a_{n+n_{0}}$, $b^{(n_0)}_{n}=  b_{n+n_{0}}$ has   eigenvalue $z$.
  Put
   \begin{equation} 
 G_{n} (z)=\sum_{m=n_{0}}^n (a_{m-1} f_{m-1}(z) f_{m}(z))^{-1},\q n\in{\Bbb Z}_{+}.
\label{eq:GE}\end{equation}

\begin{theorem}\label{GE}
Let   $z\in{\Bbb C}\setminus [-1,1]$. 
  Under assumptions \e{eq:comp} and \e{eq:LR} the sequence $g_{n} (z)$ defined by 
        \[
g_{n} (z)= f_{n}(z) G_{n} (z) 
\]
 satisfies equation \e{eq:Jy} and
     \begin{equation}
\lim_{n\to\infty} q_{n}(z)  g_{n}(z)=\frac{1}{\sqrt{z^2-1} }.
\label{eq:gas}\end{equation}
 \end{theorem}

\begin{pf}
First, we check equation \e{eq:Jy} for $g_{n}$. According to definition \e{eq:GE}, we have
 \begin{multline*}
a_{n-1} f_{n-1} G_{n-1}+ (b_{n} -z) f_{n} G_{n} +a_{n} f_{n+1}G_{n+1}
\\
=\big(a_{n-1} f_{n-1}  + (b_{n} -z) f_{n}  +a_{n} f_{n+1}\big) G_{n}
+ a_{n-1} f_{n-1} ( G_{n-1}- G_{n}) + a_{n} f_{n+1}( G_{n+1}- G_{n}) .
\end{multline*}
The first term here is zero because    equation \e{eq:Jy} is true for the  sequence $f_{n}$ Since
\[
G_{n+1}= G_n+ (a_{n} f_{n} f_{n+1})^{-1},
\]
the second and third terms equal $-f_{n}^{-1}$ and $f_{n}^{-1}$, respectively. 

Next, we prove asymptotics \e{eq:gas}. 
Recall that   $f_{n}= q_{n}u_{n} $ where $q_{n}$ is defined by formulas \e{eq:aa1},  \e{eq:re1}  and 
$u_{n}$ is constructed in Theorem~\ref{eik2}. Let us
set
\begin{equation}
 v_{m}=\z_{m-1}\z_{m} (1-\z_{m-1}\z_{m})^{-1}(a_{m-1} u_{m-1} u_{m})^{-1} ,\q t_{m}=(q_{m-1} q_{m})^{-1} .
\label{eq:hv}\end{equation}
Then 
\[
t_{m}' =  \big((\z_{m-1}\z_{m})^{-1}  -1\big) t_{m}
\]
and
\[
(a_{m-1} f_{m-1} f_{m})^{-1}=  (a_{m-1} u_{m-1} u_{m})^{-1}t_{m}=  v_{m} t'_{m}.
\]
 Integrating  by parts (see formula \e{eq:Abel}) in \e{eq:GE}, we  find that
 \begin{equation}
 G_{n} =\sum_{m=n_{0}}^n v_{m}  t_{m}' = v_{n}  t_{n+1}    -  v_{n_{0}-1}   t_{n_{0}} 
-   \sum_{m=n_{0}}^n v_{m-1}'   t_{m}.
\label{eq:GE1}\end{equation}

Let us multiply this expression by $q_{n}  f_{n} = q_{n}^2 u_{n}$ and pass to the limit $n\to\infty$. Since $v_n\to 2 \z^{2} (1-\z^{2})^{-1}$,
we see that
\begin{equation}
q_{n}^2 u_{n} v_{n}  t_{n+1} =\z_{n}^{-1}u_{n}v_{n}\to 2 \z  (1-\z^{2})^{-1}=\frac{1}{\sqrt{z^2-1} }
\label{eq:GX}\end{equation}
as $n\to \infty$.   
 It follows from  \e{eq:GXx} that $q_{n}^2 u_{n} v_{n_{0}-1}   t_{n_{0}} = O (r^{2n})$.
 The contribution
    \[
q_{n}^2 u_{n} \sum_{m=n_{0}}^n v_{m-1}'  t_{m}=  u_{n}  \sum_{m=n_{0}}^n v_{m-1}'  \frac{ q_{n }^2 }{q_{m-1}q_{m}}
\]
of the third term  in the right-hand side \e{eq:GE1} can  
be estimated by
     \begin{equation}
 \sum_{m=n_{0}}^n r^{2(n-m)} |v_{m-1}'  |\leq  r^{n} \sum_{ n_{0}\leq m<[n/2]} |v_{m-1}' |+ \sum_{ [n/2]\leq m \leq n}  |v_{m-1}' |.
\label{eq:hv1}\end{equation}
Recall that $\{ u_{n}' \}\in \ell^1 ({\Bbb Z}_{+})$ according to Remark~\ref{u'} and
$\{ \z_{n}' \}\in \ell^1 ({\Bbb Z}_{+})$ according to   definitions \e{eq:aa} and \e{eq:aa1}. It follows from \e{eq:hv} that $\{ v_{n}' \}\in \ell^1 ({\Bbb Z}_{+})$, and therefore   expression \e{eq:hv1} tends to zero as $n\to\infty$. Hence  \e{eq:GX} implies \e{eq:gas}.
 \end{pf}

 By definition \e{eq:GE}, the Wronskian \e{eq:Wr}  of $f (z)=\{f_{n}(z)\}$ and $g (z)=\{f_{n}(z) G_{n}(z)\}$ equals
  \[
\{ f(z),g(z)\}= a_{n}f_{n}(z)f_{n+1}(z)( G_{n+1}(z)- G_{n}(z))=1,
\]
whence solutions $f(z)$ and $g (z)$ are linearly independent. It follows that 
 \[
P_{n} (z)= d_{+}(z)f_{n} (z)+d_{-}(z)g_{n} (z)
\]
where 
\[
d_{+}(z)= \{ P(z),g(z)\}\q \mbox{and}  \q d_{-}(z)=- \{ P(z),f(z)\}= -\Omega(z)
\]
according to \e{eq:wei}.   Obviously,  
$d_{+}(z)\neq 0$ if $d_{-}(z)= 0$. Therefore Theorems~\ref{EIK}   and \ref{GE} imply the following result.

\begin{theorem}\label{GE1}
  Under assumptions \e{eq:comp} and \e{eq:LR}  the relation \e{eq:HSq5}
is true for all  $z\in{\Bbb C}\setminus [-1,1]$  
with convergence uniform on compact subsets of $z\in{\Bbb C}\setminus [-1,1]$.
Moreover, if $\Omega(z)=0$, then
 \begin{equation}
\lim_{n\to\infty}  q_{n}(z)^{-1} P_{n}(z)=\{ P(z),g(z)\} \neq 0.
\label{eq:GEGEx}\end{equation}
 \end{theorem}

Note that Theorem~\ref{REGB} does not follow from Theorem~\ref{GE1} because asymptotics \e{eq:HSq5} is not uniform as $z$ approaches the cut along $(-1,1)$.
 
The existence of the limit in \e{eq:HSq5} is the  classical result of the Szeg\H{o} theory. It is stated as Theorem~12.1.2 in the book  \cite{Sz} where the assumptions are imposed   on the measure $d\rho(\lambda)$; in particular, it is assumed that
$\supp\rho\subset [-1,1]$. Under short-range assumption \e{eq:Tr} asymptotic relations \e{eq:HSq5}, \e{eq:GEGEx} were established in   \cite{KS} and, by a different method, in  \cite{Y/LD}.

\subsection{Asymptotics on the continuous spectrum} 

To find asymptotic behavior of the polynomials $P_{n} (\lambda)$ for $\lambda\in (-1,1)$,  that is, on the continuous spectrum of the Jacobi  operator $J$,
we have to consider two (normalized) Jost solutions $\mathbf{f} (\lambda\pm i0)=\{ \mathbf{f}_{n} (\lambda\pm i0)\}_{n=-1}^\infty$ for
$ \lambda=\cos\theta\in (-1,1)$. Of course these two solutions are complex conjugate to each other.  
 Calculating the Wronskian \e{eq:Wr} of $\mathbf{f} (\lambda+ i0)$  and $\mathbf{f}  (\lambda- i0)$  for $n\to\infty$
 and using \e{eq:Jo1}, \e{eq:Jost1}, \e{eq:JostX},  we see that
 \begin{multline*}
\{\mathbf{f}  (\lambda+i0),  \mathbf{f} (\lambda-i0)\}= a_{n}\big(e^{-i \varphi_{n}(\lambda)}  e^{i \varphi_{n+1}(\lambda)}-e^{-i \varphi_{n+1}(\lambda)}e^{i \varphi_{n}(\lambda)}\big)+ o(1)
\\
= i   \sin \theta_{n}(\lambda)+ o(1)
= i   \sin \theta (\lambda)     =i   \sqrt{1-\lambda^2}\neq 0,
 \end{multline*}
and hence these solutions  are linearly independent.  It follows that
 \begin{equation}
P_{n} (\lambda)= \ov{c(\lambda)} \mathbf{f}_{n} (\lambda+i0)+ c(\lambda)  \mathbf{f}_{n} (\lambda-i0)
\label{eq:HH2}\end{equation}
for some complex constant $c (\lambda)$. Taking the Wronskian of this equation with $ \mathbf{f}  (\lambda+i0)$   and using notation \e{eq:nJ}, we find that 
  \begin{equation}
- c (\lambda)\{\mathbf{f}  (\lambda+i0),  \mathbf{f} (\lambda-i0)\}=\{P (\lambda),  \mathbf{f} (\lambda+ i0)\}  =  \boldsymbol{\Omega} (\lambda+i0) .
\label{eq:Hx}\end{equation} 
Thus \e{eq:HH2} leads to the same formula as   \e{eq:zF} in the short-range case.

\begin{lemma}\label{HH}
 For $ \lambda\in (-1,1)$,  the representation 
  \begin{equation}
P_{n} (\lambda)=\frac{  \boldsymbol{\Omega} (\lambda-i0)  \mathbf{f} _{n} (\lambda+i0) -  \boldsymbol{\Omega} (\lambda+i0)  \mathbf{f} _{n} (\lambda-i0)  }{i   \sqrt{1-\lambda^2}},  \q \lambda\in (-1,1), \q n=0,1,2, \ldots, 
\label{eq:HH4L}\end{equation}
 holds true.
 \end{lemma}
 
 Properties of the Wronskian \e{eq:Hx} are summarized in the following statement.
 
 \begin{theorem}\label{HX}
 The Wronskians $\boldsymbol{\Omega} (\lambda + i0)$ and $\boldsymbol{\Omega} (\lambda - i0)=
 \ov{\boldsymbol{\Omega} (\lambda + i0) }$  are continuous functions of $\lambda\in (-1,1)$ and
 \begin{equation}
\boldsymbol{\Omega} (\lambda\pm i0) \neq 0 ,\q \lambda\in (-1,1)  .
\label{eq:HH5}\end{equation}
 \end{theorem}
 
  \begin{pf}
  The functions 
 $\boldsymbol{\Omega} (\lambda\pm i0)$ are    continuous    by Corollary~\ref{JOST}.
  If $\boldsymbol{\Omega} (\lambda\pm i0)=0$, then according to \e{eq:HH4L}
    $P_{n} (\lambda)=0$ for   all $n\in{\Bbb Z}_{+}$.  However, 
 $P_0 (\lambda)=1$ for all $\lambda$.
    \end{pf}

Let us set
 \begin{equation}
\varkappa (\lambda) =2| \boldsymbol{\Omega} (\lambda+i0)|, \q
-2 \boldsymbol{\Omega} (\lambda\pm i0) =   \varkappa (\lambda) e^{\pm i \eta (\lambda) },\q \varkappa (\lambda) >0.
\label{eq:AP}\end{equation}
In the theory of short-range perturbations of the Schr\"odinger operator, the functions $\varkappa (\lambda) $ and $\eta (\lambda)$ are known as the limit amplitude and the limit phase, respectively; the function   $\eta (\lambda)$ is also called the scattering  phase or the   phase shift.    Definition \e{eq:AP}     fixes $\eta (\lambda)$  only up to a term $2\pi k$ where $k\in{\Bbb Z}$.

Combined together  relations \e{eq:Jost1} and \e{eq:HH4L} yield the   asymptotics of Bernstein-Szeg\H{o} type for the  polynomials  $P_{n} (\lambda)$. Recall that $\varepsilon_{n} $ are defined by \e{eq:eps}.

 \begin{theorem}\label{Sz}
 Let assumptions~\e{eq:comp} and \e{eq:LR} be satisfied,   let   $\lambda \in (-1,1)$ and let the phase $\varphi_{n}  (\lambda)$ be defined by formulas \e{eq:Joz} -- \e{eq:Jo1}. Then the  polynomials $P_{n} (\lambda)$ have asymptotics 
   \begin{equation}
 P_{n} (\lambda)=  \varkappa ( \lambda) (1-\lambda^2)^{-1/2} \sin (\varphi_{n}  (\lambda) + \eta(\lambda) ) + O(\varepsilon_{n})  
\label{eq:Sz}\end{equation}
as $n\to\infty$. Relation \e{eq:Sz} is uniform in $\lambda$ on compact subintervals of $(-1,1)$.
 \end{theorem}
 
 The phase shifts $\eta (\lambda)$  in \e{eq:Sz} and $\pi\xi(\lambda)$ in the short-range formula \e{eq:zF} play the same roles.  They depend on the precise values of the coefficients $a_{n}$ and $b_{n}$ for all $n$ and hence   cannot be found from their asymptotic behavior as $n\to\infty$. Under additional assumptions the growing part $\varphi_{n}  (\lambda) $ of the phase in \e{eq:Sz} can be made more explicit.

 \subsection{Hilbert-Schmidt perturbations}
 
 In addition to   \e{eq:LR}, assume now that condition 
 \begin{equation}
 \sum_{n=0}^\infty  (\alpha_{n}^2 + b_{n}^2 )<\infty,
\label{eq:HS}\end{equation}
 is satisfied,
 that is, $V=J-J_{0}$ is a Hilbert-Schmidt operator. Then  asymptotic formulas  of Theorems~\ref{GE1} and \ref{Sz} can be made more explicit. We proceed from the following elementary assertion.
 
 \begin{lemma}\label{HSq}
 Let $z\neq \pm 1 $. 
 Under assumption \e{eq:HS} there exists a finite limit
  \begin{equation}
\lim_{n\to\infty}  \Big(\z(z)^{-n} \exp \big(-\frac{1}{\sqrt{z^2-1}}\sum_{m=0}^{n-1}(2 z\alpha_{m} + b_m)  \big)q_{n}(z)\Big)\neq 0.
\label{eq:HSq}\end{equation}
\end{lemma}

\begin{pf}
It follows (cf. \e{eq:aa1-}) from \e{eq:aa}, \e{eq:aa1} that
 \begin{equation}
\z_{m}-\z   = (z_{m}-z)
 \Big(1-\frac{z_{m}+z}{\sqrt{z_{m}^2-1}+\sqrt{z^2-1}}\Big).
\label{eq:HSq1}\end{equation}
Since
\[
z_{m}-z=-2 z\alpha_{m}-b_{m}+O(\alpha_m^2+b^2_m)
\]
and
\[
1-\frac{z_{m}+z}{\sqrt{z_{m}^2-1}+\sqrt{z^2-1}}=-\frac{\z}{\sqrt{z^2-1}}+O(\sqrt{\alpha_m^2+b^2_m}),
\]
equality \e{eq:HSq1} implies that
 \[
\frac{\z_{m}}{\z}   = 1+ \frac{2 z\alpha_{m}+b_{m}}{\sqrt{z^2-1}} +O(\alpha_m^2+b^2_m)
=\exp \Big(\frac{2 z\alpha_{m}+b_{m}}{\sqrt{z^2-1}}\Big)\big(1+O(\alpha_m^2+b^2_m)\big).
\]
Taking the product over $m=0,1,\ldots, n -1$ and using condition  \e{eq:HS}, we arrive at  \e{eq:HSq}.
\end{pf}

Now the following statement is a direct consequence of Theorem~\ref{GE1}.

 \begin{theorem}\label{SzHSq}
  Let assumptions     \e{eq:LR} and  \e{eq:HS} be satisfied, and  let $z\in  {\Bbb C}\setminus [-1,1] $.  Then there exist  finite limits
  \begin{equation}
\lim_{n\to\infty}  \Big(\z(z)^{n} \exp \big(\frac{1}{\sqrt{z^2-1}}\sum_{m=0}^{n-1}(2 z\alpha_{m} + b_m)  \big)P_{n}(z)\Big)\neq 0
\label{eq:HSq3}\end{equation}
if $z$ is not an eigenvalue of the  operator $J$ and
 \begin{equation}
\lim_{n\to\infty}  \Big(\z(z)^{-n} \exp \big(-\frac{1}{\sqrt{z^2-1}}\sum_{m=0}^{n-1}(2 z\alpha_{m} + b_m)  \big)P_{n}(z)\Big)\neq 0
\label{eq:HSp3}\end{equation}
if $z$ is   an eigenvalue of   $J$.
 \end{theorem}
 
  \begin{corollary}\label{SzHSq1}
  Suppose additionally that the conditions
    \begin{equation}
\sum_{n=0}^\infty \alpha_{n}< \infty \q \mbox{and} \q\sum_{n=0}^\infty b_{n}< \infty 
\label{eq:HSq4}\end{equation} 
$($these series should be convergent but perhaps not absolutely$)$ are satisfied. Then the exponential factors in  \e{eq:HSq3} and  \e{eq:HSp3} may be omitted, that is, the limits in \e{eq:HSq5} and \e{eq:GEGEx} exist.
\end{corollary}

 
 We do not know whether 
relations \e{eq:HSq3} and \e{eq:HSp3} remain true under the only assumption   \e{eq:HS}.

  In some cases the exponential factors in  \e{eq:HSq3} and  \e{eq:HSp3} can be  simplified.
    
\begin{example}\label{HSe}
Let  conditions \e{eq:LR1} be satisfied with some   $r_{1}, r_{2}\in (1/2,1)$. 
Then 
  \begin{multline}
\sum_{m=0}^n (2z\, \alpha_{m} + b_m)=2 z\alpha (1-r_{1})^{-1}n^{1- r_{1}}+ b
(1-r_{2})^{-1}n^{1- r_{2}}
\\
+ 2 z  \alpha\boldsymbol{\gamma}_{r_{1}}+
b \boldsymbol{\gamma}_{r_{2}} + \sum_{m=0}^\infty (2 z  \tilde{\alpha}_m+ \tilde{b}_m)+  o(1)
\label{eq:HSp}\end{multline}
where $\boldsymbol{\gamma}_r-(1-r)^{-1}$  is the Euler-Mascheroni constant.
With a natural modification, expression \e{eq:HSp} remains true if $r_{j} =1$ for one or both $j$. In this case
$(1-r_{j})^{-1}n^{1-r_{j}}$ should be replaced by $\ln n$ and  $\boldsymbol{\gamma}_1$ is    the Euler-Mascheroni constant.
  \end{example}

Let us now discuss relation \e{eq:Sz}. Similarly to Lemma~\ref{HSq}, we have

\begin{lemma}\label{HS}
 Under assumption \e{eq:HS} there exists a finite limit
  \begin{equation}
\lim_{n\to\infty}  \big(\varphi_n (\lambda) - n\theta-(\sin\theta)^{-1}\sum_{m=0}^{n-1} (2\cos \theta\, \alpha_{m} + b_m)\big)=: \gamma (\lambda).
\label{eq:HS1}\end{equation}
\end{lemma}

\begin{pf}
It follows from \e{eq:Joz}, \e{eq:Joa} that
 \[
\theta_m=\arccos\frac{\lambda-b_m}{2a_m}=\theta + (\sin\theta)^{-1} (2\cos \theta\, \alpha_m + b_m)+O(\alpha_m^2+b^2_m).
\]
Taking the sum over $m=0,1,\ldots, n -1$ and using condition  \e{eq:HS}, we arrive at  \e{eq:HS1}.
\end{pf}

Now the following statement is a direct consequence of Theorem~\ref{Sz}.

 \begin{theorem}\label{SzHS}
  Let assumptions     \e{eq:LR} and  \e{eq:HS} be satisfied. Then
for  $\lambda \in (-1,1)$, the asymptotic formula 
   \begin{equation}
 P_{n} (\lambda)=   \varkappa ( \lambda) (1-\lambda^2)^{-1/2}  \sin \big(n\theta + (\sin\theta)^{-1}\sum_{m=0}^{n-1} (2\cos \theta\, \alpha_{m} + b_m) +\gamma (\lambda)+ \eta(\lambda)\big) + o(1)  
\label{eq:SzHS}\end{equation}
holds as $n\to\infty$. Relation \e{eq:SzHS} is uniform in $\lambda$ on compact subintervals of $(-1,1)$.
 \end{theorem}
 
   Under assumption \e{eq:LR1} the phase in \e{eq:SzHS} can be simplified if one takes relation \e{eq:HSp} (where $z$ is replaced by $\cos\theta$) into account.
  
  Of course   formulas \e{eq:HSq3} and \e{eq:SzHS} are consistent with asymptotic  formulas   for Pollaczek
  polynomials in the Appendix in the book \cite{Sz}.

Without any additional assumptions, Hilbert-Schmidt perturbations $V$ of the operator $J_{0}$
were investigated in the deep papers \cite{KS} and \cite{DS}.  In   \cite{KS}, necessary and sufficient conditions in terms of the spectral measure $d\rho(\lambda)$ of the operator $J=J_{0} + V$ were found for $V$ to be in the   Hilbert-Schmidt class.  Asymptotic behavior of the corresponding polynomials $P_{n} (z)$  was studied in \cite{DS}. It was proved in   Theorem~5.1 that
  the   limit of $ \z(z)^{n} P_{n}  (z)$ as $n\to\infty$ exists 
if and only if  conditions \e{eq:HS} and \e{eq:HSq4} are satisfied.
As shown in   Theorem~8.1 of \cite{DS},    assumptions \e{eq:HS},  \e{eq:HSq4} are sufficient also for the validity of  formula  \e{eq:zC} but only in some {\it averaged}  sense. Such a regularization seems to be necessary since under  these  assumptions    
 the structure of the essential spectrum of the operator $J$ can be quite wild.  

 Condition \e{eq:LR} accepted in this paper is different in nature from \e{eq:HS},  \e{eq:HSq4}. On the one hand, it excludes too strong oscillations of the coefficients $\alpha_{n}$, $b_{n}$ but, on the other hand, it permits their arbitrary slow decay as $n\to\infty$.

\section{Spectral theory of Jacobi operators}

Here we show that the spectrum of the Jacobi operator $J$ on the interval $(-1,1)$ is absolutely continuous and the corresponding weight $w(\lambda)$ is expressed via the Jost function   by the formula \e{eq:wei1}.
It follows that $w(\lambda)$ is
a continuous strictly positive function of $\lambda\in (-1,1)$.

\subsection{Resolvent}
 
  First, we construct
the resolvent $R(z)=(J -z )^{-1}$ of the operator $J$. Recall that   $P  (z)= \{P_{n}(z)\} $ and $f(z)=\{f_{n} (z)\}$ are, respectively,  the polynomial and the Jost solutions of the   equation \e{eq:Jy} and 
$\Omega(z)$ is their Wronskian \e{eq:wei}. We denote by $e_{n}$, $n\in {\Bbb Z}_{+}$, the canonical basis in $\ell^2 ({\Bbb Z}_{+})$.   
The proof below is almost the same as in the short-range case (cf. Lemma~2.6 in \cite{Y/LD}).

 \begin{lemma}\label{res}
 For all $n,m\in{\Bbb Z}_{+}$, we have
  \begin{equation}
(R (z)e_{n}, e_{m})= \Omega(z)^{-1} P_{n} (z) f_{m}(z),\q \Im z \neq 0, 
\label{eq:RR}\end{equation}
if $n\leq m$ and $(R (z)e_{n}, e_{m})=(R (z)e_m, e_n)$.
 \end{lemma} 
 
  \begin{pf} 
  We will show that
   the operator $R(z)$   defined by relation \e{eq:RR} is the resolvent of $J$. We have
    \begin{equation}
  \Omega(z)( R (z)u)_{n} =  f_{n}(z) A_{n}(z)+   P_{n}(z) B_{n}(z)
\label{eq:RR1}\end{equation}
where
   \begin{equation}
A_{n}(z) =\sum_{m=0}^n  P_{m} (z) u_{m},\q    B_{n}(z) =\sum_{m=n+1}^\infty  f_{m} (z) u_{m},
\label{eq:RR2}\end{equation}
at least for all sequences $u=\{u_{n}\}$ with a  finite number of non-zero components $u_{n}$. In this case $R (z)u\in \ell^2 ({\Bbb Z}_{+})$ because $f_{n} (z)\in \ell^2 ({\Bbb Z}_{+})$ if $\Im z\neq 0$. 

Our goal is to check that $(J-z)R (z)u=u$. 
It follows from definition \e{eq:ZP+}  of  the Jacobi operator $J$ and formula \e{eq:RR1} that
   \begin{multline}
  \Omega ( (J-z) R  u)_{n} =  a_{n-1}\big(f_{n-1}A_{n-1}+P_{n-1} B_{n-1} \big)
  \\
  + (b_{n}-z) \big(f_{n}A_{n}+ P_{n} B_{n}  \big)+ a_{n}\big(f_{n+1}A_{n+1} +P_{n+1} B_{n+1}    \big).
\label{eq:RR3}
\end{multline}
According to \e{eq:RR2} we have
\[
f_{n-1}A_{n-1}+P_{n-1} B_{n-1}=f_{n-1}(A_{n}-P_{n} u_{n})+  P_{n-1}(B_{n}+f_{n} u_{n})
\]
and
\[
f_{n+1} A_{n+1} +P_{n+1} B_{n+1} =
f_{n+1}A_{n}  +  P_{n+1}B_{n} .
\]
Let us substitute these expressions into the right-hand side of \e{eq:RR3} and observe that
the coefficients at $A_{n}$ and $B_{n}$ equal zero by virtue of equation \e{eq:Jy} for $\{f_{n}\}$ and $\{P_{n}\}$, respectively. It  follows that
\[
 ( (J-z) R  u)_{n} = \Omega^{-1}  a_{n-1}  (- P_{n}  f_{n-1} +    f_{n} P_{n-1}) u_{n}  =  u_{n} 
\]
whence  $R(z)=(J-z)^{-1}$. In particular, the operator $R(z) $ defined by \e{eq:RR} is bounded in the space $  \ell^2 ({\Bbb Z}_{+})$.
    \end{pf}
  
 In view of Theorem~\ref{EIK},   $f_{n} (z)$, $n=-1, 0,1,\ldots$,  and $\Omega (z)$ are analytic    functions of $z\in{\Bbb C}\setminus {\Bbb R}$, and they are continuous up to the cut along ${\Bbb R}$ with possible exception of the points $z=\pm 1$. The values $f_{m} (\lambda\pm i0)$ and $\Omega (\lambda\pm i0)$ on the upper and lower edges of the cuts along $(-\infty,- 1)$ and $(1,\infty)$ are, in general, different. Note however the following obvious
 
 \begin{remark}\label{AR}
If $\lambda\in (-\infty,- 1)\cup (1,\infty)$ and $\Omega (\lambda\pm i0)\neq 0$, then
 \begin{equation}
\frac{f_{m} (\lambda + i0)}{\Omega(\lambda + i0)}= \frac{f_{m} (\lambda - i0)}{\Omega (\lambda - i0)}.
\label{eq:AR}\end{equation}
Indeed, consider the sequence
\[
\Delta_{m} (\lambda)=f_{m} (\lambda + i0)\Omega(\lambda - i0)- f_{m} (\lambda - i0)\Omega(\lambda + i0).
\]
It satisfies equation \e{eq:Jy} where $z=\lambda$, belongs to $\ell^2 ({\Bbb Z}_{+})$ because $\lambda\not\in [-1,1]$ and $\Delta_{-1} (\lambda)=0$ by definition \e{eq:wei} of $\Omega (\lambda\pm i0)$. Since $\lambda$ is not an eigenvalues of $J$, we see that $\Delta_{m} (\lambda)=0$ for all $m$ which is equivalent to \e{eq:AR}.
It follows from \e{eq:AR} that the function $f_{m}  (z)/ \Omega (z)$ is analytic in ${\Bbb C}  \setminus [-1,1]$ with poles at eigenvalues of the operator $J$. This is consistent with formula \e{eq:RR} since the matrix elements $ (R(z) e_{n}, e_{m})$ are analytic functions of  $z\in {\Bbb C}  \setminus [-1,1]$.
 \end{remark}

 Taking also \e{eq:HH5} into account, we obtain the following result.

 \begin{theorem}\label{AC}
  Let assumptions \e{eq:comp} and \e{eq:LR}  hold. Then
  \begin{enumerate}[{\rm(i)}]
 \item
The  resolvent $R(z)=(J-z)^{-1}$ of the Jacobi operator $J$ is an integral operator with matrix elements  \e{eq:RR}.
 For all $n,m\in{\Bbb Z}_{+}$,  it    is an analytic function of $z\in {\Bbb C}\setminus [-1, 1]$ with simple poles at eigenvalues of the operator $J$. A point
$z\in {\Bbb C}\setminus [-1, 1]$ is an eigenvalue of $J$ if and only if $\Omega(z)=0$.

 \item
 For all $n,m\in {\Bbb Z}_{+}$, the functions $(R(z) e_{n}, e_{m})$ are  continuous in $z$ up to the cut $[-1,1]$ except, possibly, the points $\pm 1$.
 
  \item
Estimates
\[
|(R (z)e_{n}, e_{m})|\leq C | \Omega(z)|^{-1} \big|q_{m}(z)/q_{n}(z)\big|\leq C_{1}<\infty, \q n\leq m,
\]
are true  with some positive constants that do not depend on $n$, $m$ and on $z$ in compact subsets of the  set $\clos(\Bbb C\setminus  [-1,1])$ as long as they are away from the points $\pm 1$.
 \end{enumerate}
 \end{theorem}
 
  The statement (ii)   is known as the limiting absorption principle. It implies
 
 \begin{corollary}\label{Rc}
The spectrum of the operator $J$ on the interval $(-1,1)$ is absolutely continuous.
 \end{corollary}

Let us now consider the spectral projector $E(\lambda)$ of the operator $J$. By the Cauchy-Stieltjes formula, its matrix elements satisfy the identity
 \begin{equation}
2\pi i \frac{d(E (\lambda)e_n, e_m)} {d\lambda}=  (R (\lambda+ i0)e_n, e_m)-(R (\lambda- i0)e_n, e_m).
\label{eq:Priv}\end{equation}
The following assertion is a direct consequence of  Theorem~\ref{AC}, part (ii).

\begin{corollary}\label{RE}
 For all $n,m\in {\Bbb Z}_{+}$, the functions $(E(\lambda) e_{n}, e_{m})$ are  continuously differentiable in $\lambda\in (-1,1)$.
   \end{corollary}

 We emphasize that the points $1$ and $-1$ may be eigenvalues of $J$; see Example~4.15 in \cite{Y/LD}.
 
 Theorem~\ref{AC} can also be obtained by the Mourre method \cite{Mo1}. It was applied to Jacobi operators in \cite{BdM}; to be precise, the problem in the space $\ell^2 ({\Bbb Z})$ was considered in \cite{BdM}, but this is of no importance. However,  the Mourre method does not exclude eigenvalues of $J$ embedded in its continuous spectrum although it shows that  these eigenvalues   may accumulate to the points $1$ and $-1$ only. Note also that very general conditions of the absolute continuity of spectrum were obtained in \cite{Stolz} by the subordinacy method of \cite{GD}.

 \subsection{Spectral measure}
 
Now we are in a  position to calculate  the spectral family  $d E(\lambda)$ of the   operator $J $.
Let us proceed from the identity \e{eq:Priv}.
 Using notation \e{eq:nJ} and \e{eq:Hx} we can rewrite
  formula \e{eq:RR} as
\begin{equation}
(R (\lambda\pm i0)e_{n}, e_{m})=  \boldsymbol{\Omega} (\lambda\pm i 0)^{-1} P_{n} (\lambda) \mathbf{f}_{m}(\lambda\pm i 0),\q n\leq m,
\label{eq:RRpm}\end{equation}
where
\[
\mathbf{f}_{m}(\lambda - i 0)=\ov{\mathbf{f}_{m}(\lambda + i 0)}\q \mbox{and} \q
\boldsymbol{\Omega}(\lambda - i 0)=\ov{\boldsymbol{\Omega}(\lambda + i 0)} .
\]
Substituting  expression \e{eq:RRpm} into \e{eq:Priv}, we find that
\[
2\pi i \frac{d(E (\lambda)e_n, e_m)} {d\lambda}= P_{n} (\lambda) \frac{ \boldsymbol{\Omega} (\lambda-i 0)  \mathbf{f}_{m}(\lambda+i 0)    - \boldsymbol{\Omega} (\lambda+i 0) \mathbf{f}_{m}(\lambda-i 0)    }{| \boldsymbol{\Omega} (\lambda+i 0) |^2}.
\]

Combining this representation with  formula \e{eq:HH4L} for $P_m (\lambda) $, we obtain    the following result.

 \begin{theorem}\label{SF}
  Let assumptions \e{eq:comp} and \e{eq:LR}  hold. Then, for all $n,m\in{\Bbb Z}_{+}$ and $\lambda\in(-1,1)$, we have the representation
  \[
 \frac{d(E (\lambda)e_n, e_m)} {d\lambda}= (2\pi)^{-1}\sqrt{1-\lambda^2}  |\boldsymbol{\Omega} (\lambda\pm i0) |^{-2} P_{n} (\lambda) P_m (\lambda)   .
\]
 In particular,
the spectral measure of the operator $J$ equals
  \begin{equation}
d\rho(\lambda):=d(E (\lambda)e_0, e_0)= w(\lambda) d\lambda,\q \lambda\in(-1,1),
\label{eq:SFx}\end{equation}
where the weight  $w(\lambda)$ is given by the formula 
 \begin{equation}
w(\lambda)=  (2\pi)^{-1}\sqrt{1-\lambda^2}  \, | \boldsymbol{\Omega} ( \lambda\pm i0) |^{-2}     
\label{eq:SF1}\end{equation}
$($the right-hand sides here do not depend on the sign$)$.
 \end{theorem} 
 
 According to \e{eq:SF1}  the amplitude factors in \e{eq:zF} and \e{eq:Sz} are the same. 

 Putting together Theorem~\ref{HX} and formula \e{eq:SF1}, we arrive at the next result.
 
  \begin{theorem}\label{SFr}
  Under assumptions \e{eq:comp} and \e{eq:LR} the weight $w (\lambda)$ is a continuous strictly positive function of $\lambda\in (-1,1)$.
 \end{theorem}

 Note that this result was earlier obtained in \cite{M-N} by specific methods of the orthogonal polynomials theory.

According to \e{eq:WRfr} for the operator $J_{0}$, we have
 \[
 \boldsymbol{\Omega}_{0}(\lambda\pm i0)=\Omega_{0}(\lambda\pm i0)=-2^{-1} (\lambda\pm i\sqrt{1-\lambda^2}),
 \]
 and hence expressions  \e{eq:SFx}, \e{eq:SF1}  reduce to \e{eq:fr}.

  Observe that Theorem~\ref{EIK}  does not give any information on the behavior of the Jost function $\Omega(z)$ as $z\to\pm 1$. However relation \e{eq:SF1} implies that the normalized function \e{eq:nJ} obeys an estimate
 \[
  \int_{-1}^1\sqrt{1-\lambda^2}   | \boldsymbol{\Omega} ( \lambda\pm i0) |^{-2} d\lambda
  = 2\pi 
  \int_{-1}^1 w(\lambda) d\lambda\leq 2\pi ,
\]
and hence $ \boldsymbol{\Omega} ( \lambda\pm i0) $ cannot vanish too rapidly as $\lambda\to 1-0$ and $\lambda\to -1+0$ (even if $1$ or $-1$ are eigenvalues of $J$).
  For example, the behavior $\boldsymbol{\Omega} ( \lambda+i0) \sim c_{\pm}  (\lambda\mp 1)$ with $c_{\pm}  \neq 0$ is excluded.  

In view of \e{eq:AP}, \e{eq:SF1} the amplitude in  \e{eq:Sz}   can  be written as
    \begin{equation}
 \varkappa (\lambda)( 1-\lambda^2)^{-1/2}= (2 /\pi)^{1/2}  (1-\lambda^2)^{-1/4} w(\lambda)^{-1/2} . 
\label{eq:Sz1}\end{equation}
Substituting this expression into \e{eq:Sz}, we can reformulate Theorem~\ref{Sz}
in a form more common for the orthogonal polynomials literature.

 \begin{theorem}\label{SW}
 Let assumptions \e{eq:comp} and \e{eq:LR} be satisfied,   let   $\lambda \in (-1,1)$ and let the phases $\varphi_{n}  (\lambda)$ and $\eta(\lambda)$ be defined by formulas \e{eq:Jo1} and \e{eq:AP}, respectively. Then the  polynomials $P_{n} (\lambda)$ have asymptotics 
   \begin{equation}
 P_{n} (\lambda)=   (2/\pi)^{1/2}  (1-\lambda^2)^{-1/4} w(\lambda)^{-1/2}  \sin (\varphi_{n}  (\lambda) + \eta(\lambda) ) + O( \varepsilon_{n})  
\label{eq:SW}\end{equation}
as $n\to\infty$. Relation \e{eq:SW} is uniform in $\lambda$ on compact subintervals of $(-1,1)$.
 \end{theorem}
 
 Formula \e{eq:SW} is a generalization of \e{eq:zF} and reduces to it if the short-range condition \e{eq:Tr} is satisfied.
 
 Of course we can also replace the amplitude factor in \e{eq:SzHS} by its expression \e{eq:Sz1}.

 \subsection{Singular weights}

  As was mentioned in Sect.~1.2, results on orthogonal polynomials $P_{n}  (z)$ can be deduced from some conditions on either the recurrent coefficients $a_{n} , b_{n}$ or on the corresponding spectral measure $d\rho(\lambda)$. We proceeded from conditions on $a_{n} , b_{n}$, but our approach gives automatically also some results about the regularity of $d\rho(\lambda)=w(\lambda) d\lambda$ (see Theorem~\ref{SFr}).
    
  Weights  $w (\lambda)$ with  singularities or zeros inside $(-1,1)$   change both the asymptotic behavior of the recurrent coefficients $a_{n} , b_{n}$ and of the orthogonal polynomials $P_{n}  (z)$. For example, even weights $w (\lambda)$ behaving like $\kappa |\lambda|^\gamma$ where $\gamma > -1$, $\kappa>0$ as $\lambda\to 0$ were  
      investigated in the paper \cite{Nev}. Such weights are either singular at the point $\lambda=0$ if $\gamma <0$ or $w(0)=0$ if $\gamma >0$. It was shown in \cite{Nev}, Theorem~7.4,  that the corresponding Jacobi coefficients $a_{n}$ satisfy the asymptotic relation
   \[
   a_{n} =1/2 + (-1)^n \gamma  / (4n) + o(1/n)
   \]
   (the coefficients $b_{n}=0$ if the weight $w(\lambda)$ is even).
   Since $|a_{n}'|\sim |\gamma|/ (2n)$, the condition \e{eq:LR} is now violated.   For such weights, the asymptotic behavior of the polynomials $P_{n}  (\lambda)$ in a neighborhood of the point $\lambda=0$ differs from \e{eq:LRz} and  from \e{eq:zF}, in particular.
   The results of \cite{M-F}  for weights with a jump singularity are morally similar to those in \cite{Nev}. 
    Thus, the long-range condition  \e{eq:LR} 
   is close to necessary even for our   results on the weight $w(\lambda)$.
   
    In our case  the weight $w(\lambda)$ is a continuous positive function. Nevertheless the classical asymptotics \e{eq:zF}  breaks down  since an additional phase shift $\Phi_{n}  (\lambda)$ appears in \e{eq:LRz}.  This is quite similar to long-range scattering for the Schr\"odinger equation.
  We emphasize, however, that   asymptotics \e{eq:LRz} obtained under assumption \e{eq:LR}  is, in some sense, more regular than  the asymptotics of $P_{n}  (\lambda)$ in \cite{Nev,   M-F}.

 

The author states that there is no conflict of interest.


\begin{thebibliography}{99}
 
 \bibitem{AKH} N.~Akhiezer,
\emph{The classical moment problem and some related questions in analysis},  Oliver and  Boyd, Edinburgh and London, 1965.


  \bibitem{Assche} W.~Van Assche, Asymptotics for orthogonal polynomials and three-term recurrences,
Orthogonal polynomials (Columbus, OH, 1989) NATO Adv. Sci. Inst. Ser. C Math.
Phys. Sci. 294 (1990), 435-462.
  
    \bibitem{Va-As} W.~Van Assche and J.~S.~Geronimo, Asymptotics of the  
 orthogonal polynomials on and off the essential spectrum, Journal   Appr. Theory {\bf 55} (1988), 220-231.

  

 \bibitem {Bern} S. Bernstein,
   Sur les polyn\^omes orthogonaux relatifs  \`a un segment fini, Journal de Math\'ematiques  {\bf 9}  (1930), 127-177; 
  {\bf 10}  (1931), 219-286.

 

\bibitem {B-L}S. Bodine and D. Lutz,
Asymptotic integration
of differential and difference
equations, LNM 2129, Springer, 2015.
 

\bibitem {BdM} A.  Boutet de Monvel and J. Sahbani,  
Anisotropic Jacobi matrices with absolutely continuous spectrum,
C. R. Acad. Sci. Paris Sér. I Math. {\bf 328}, No. 5  (1999), 443-448. 

\bibitem {Case} K.~M.~Case,  Orthogonal polynomials from the viewpoint of scattering theory, Journal Math. Phys {\bf 15}  (1974), 2166-2174. 

\bibitem {CoLe} E. A. Coddington and N. Levinson, {\it  Theory of ordinary differential equations},
McGraw-Hill, New York, 1955.
  
 
  


  \bibitem {GD} D.~Gilbert and D. B.~Pearson, On subordinacy and analysis of the spectrum of one-dimensional
   Schr\"odinger  operators, J. Math. Anal. Appl. {\bf 128}, no. 2 (1987), 30-56.
    
   
     
 
       
       

 
   
  
  
 

 



   


 
 







\bibitem {DS} D.~Damanik and B.~Simon,   Jost functions and Jost solutions for Jacobi matrices, I. A necessary and sufficient condition for Szeg\H{o} asymptotics, Invent. Math. {\bf 165}  (2006), 1-50.

\bibitem {Deift}  P.~Deift, {\it  Orthogonal polynomials and random matrices. A Riemann-Hilbert approach}, NYU lectures, AMS, 2000.

 

 \bibitem{BE} A.~Erd\'elyi,  W.~Magnus, F.~Oberhettinger, F.~G.~Tricomi, \emph{Higher transcendental functions}, Vol. 1, 2,   McGraw-Hill, New York-Toronto-London, 1953.


\bibitem {Gon} A. A.~Gon\v{c}ar, On convergence of Pad\'e approximants for some classes of meromorphic functions,   Math. USSR Sb. {\bf 26}  (1975), 555-575.


    \bibitem{Ism} M. E. H.~Ismail, {\em Classical and quantum orthogonal polynomials in one variable}, Cambridge University Press, Cambridge, 2005.
    
    \bibitem {Jost} R. Jost, {\em \"Uber die falschen Nullstellen der Eigenwerte des $S$-matrix}, Helv. Phys. Acta  {\bf 20} (1947), 250-266.
 
 
    
 

  \bibitem{KS}     R.~Killip and B.~Simon,   Sum rules for Jacobi matrices and their applications to spectral theory,  Ann.
of Math., {\bf 158} (2003), 253-321.
 
 \bibitem {LL} L. D.~ Landau and E. M.~Lifshitz, {\it Quantum mechanics}, Pergamon Press, 1965.
 
  \bibitem {Lub} D.~S.~Lubinsky,  A survey of general orthogonal polynomials for weights on finite and infinite intervals,  Acta Appl. Math.  {\bf 10},   (1987),  237-296. 
  
   \bibitem {M-N} A.~M\'at\'e and  P.~Nevai,  Orthogonal polynomials and absolutely continuous measures, In: Approximation Theory IV, 611-617 (C.~K.~Chui, L.~L.~Schumaker, J.~ D. Ward, eds.), New York: Academic Press, 1983.
   
    
 
 \bibitem {Mate} A.~M\'at\'e,  P.~Nevai, and V.~Totik,    Asymptotics for orthogonal polynomials defined by a recurrence relation,  
 Constr. Approx. {\bf 1},   (1985),  231-248. 

  \bibitem{M-F} F.~Moreno, A.~Martinez-Finkelshtein and V.~L.~Sousa,  Asymptotics of orthogonal polynomials for a weight with a jump on $[-1,1]$, Constr. Approx. {\bf 33}, No. 2,  (2011),  219-263. 

    
\bibitem {Mo1}E.~Mourre,  Absence of singular spectrum for certain self-adjoint operators, Comm.
Math. Phys. {\bf 78} (1981), 391-400. 
    
    \bibitem {Nab} S. N.~Naboko and S. I.~Yakovlev,   Discrete  Schr\"odinger   operator. 
    The point spectrum on the continuous one, Saint-Petersburg Math. Journal {\bf 4}, No. 3 (1993), 559-568.

      \bibitem {Nev} P.~Nevai,   {\em Orthogonal polynomials}, Memoirs of the AMS {\bf 18}, No. 213, Providence,
  R. I., 1979.
  
  \bibitem {Nik} E. M.~Nikishin,   Discrete Sturm-Liouville operators and some problems of function theory, J. Sov. Math. {\bf 35}  (1986), 2679-2744.

  

   \bibitem {Pos} J.~P\"oschel,   Examples of discrete  Schr\"odinger   operators with 
    pure point spectrum, Comm. Math. Phys. {\bf 88}, No. 3 (1983), 447-463.
    
       \bibitem {Stolz} G.~Stolz,   Spectral theory for slowly oscillating potentials I. Jacobi matrices, Manuscripta Math.   {\bf 84}   (1994), 245-260.
 
    \bibitem{Sz} G.~Szeg\H{o},  {\em Orthogonal polynomials}, Amer. Math. Soc., Providence, R. I., 1978.
 
 
 \bibitem {Tot}V.~Totik,  Orthogonal polynomials, Surveys in Appr. Theory {\bf 1} (2005), 70-125. 
 
   \bibitem {Weidmann} J.~Weidmann, {\it Lineare Operatoren in Hilbertr\"aumen, Teil II: Anwendungen}, Teubner Verlag, 2003. 
 
  


 


  
 
  
   \bibitem{YA} D. R. Yafaev, {\em Mathematical scattering theory: Analytic  theory}, Amer. Math. Soc.,   Providence,
  R. I., 2010.
  
  
 
  \bibitem{Y/LD} D. R. Yafaev,   Analytic scattering theory for Jacobi operators and Bernstein-Szeg\H{o}
  asymptotics of orthogonal polynomials,  Rev. Math. Phys. {\bf 30}, No. 8  (2018),  1840019.

   \bibitem{Y-LR} D. R. Yafaev,   A note on the Schr\"odinger operator with a long-range potential,      
    Letters  Math. Phys. {\bf 109}, No. 12  (2019),  2625-2648.
    
       \end{thebibliography}
  \end{document}